\input epsf
\documentclass[12pt]{article}
\usepackage{pb-diagram}
\begin{document}
\def\qq{{\bf Q}}
\def\rr{{\bf R}}
\def\zz{{\bf Z}}
\def\cc{{\bf C}}
\def\call{{\cal L}}
\def\calp{{\cal P}}
\def\calc{{\cal C}}
\def\text#1{{\mbox{#1}}}

\def\mapright#1{\smash{ \mathop{\longrightarrow}\limits^{#1}}}
\def\maptwo#1#2{\smash{
\mathop{\longrightarrow}\limits^{#1}_{#2}}}

\def\sqr#1#2{{\vcenter{\hrule height.#2pt 
\hbox{\vrule width.#2pt height #1pt \kern#1pt \vrule width.#2pt}
\hrule height.#2pt}}}
\def\square{\mathchoice\sqr{5.5}4\sqr{5.0}4\sqr{4.8}3\sqr{4.8}3}
\def\qed{\hskip4pt plus1fill\ $\square$\par\medbreak} \font\fourteenrm=cmr10
scaled\magstep1

\newtheorem{theorem}{Theorem}[section]
\newtheorem{corollary}[theorem]{Corollary} 
\newtheorem{lemma}[theorem]{Lemma}
\newtheorem{proposition}[theorem]{Proposition}

\newtheorem{xmp}[theorem]{Example}
\newtheorem{rmk}{Remark}
\newtheorem{rmks}{Remarks}
\newtheorem{dfn}[theorem]{Definition}
\newtheorem{exe}{Exercise}

\newenvironment{example}{\begin{xmp} \rm}{\end{xmp}}
\newenvironment{remark}{\begin{rmk} \rm}{\end{rmk}}
\newenvironment{remarks}{\begin{rmks} \rm}{\end{rmks}}
\newenvironment{definition}{\begin{dfn} \rm}{\end{dfn}}
\newenvironment{exercise}{\begin{exe} \rm}{\end{exe}}


\def\pictu #1 by #2 (#3 scaled #4){
      \vbox to #2{
      \hrule width #1 height 0pt depth 0pt
      \vfill
      \special{picture #3 scaled #4}}}


\def\centerpicture #1 by #2 (#3 scaled #4){
   \dimen0=#1 \dimen1=#2                                                
    \divide\dimen0 by 1000 \multiply\dimen0 by #4        
    \divide\dimen1 by 1000 \multiply\dimen1 by #4
         \noindent                                                                      
         \vbox{                                                                          
            \hspace*{\fill}                                                         
            \pictu \dimen0 by \dimen1 (#3 scaled #4)       
            \hspace*{\fill}                                                         
            \vfill}}                                                                      


\newcommand{\commandname}{\centerpicture  dim1 by dim2 (picurename scaled scalefactor)}

\begin{center}
{\bf \large THE GASSNER REPRESENTATION }
\end{center}

\begin{center}{\bf \large FOR STRING
LINKS}
\end{center}

\begin{center}{Paul Kirk, Charles Livingston, and Zhenghan Wang}
\end{center}
\begin{center}April 29, 1998\end{center}

\vskip8ex

\noindent {\bf Abstract}  The Gassner representation of the pure braid group to 
$GL_n(\zz[\zz^n])$ can be extended to give a representation of the concordance group
of $n$-strand string links to $GL_n(F)$, where $F$ is the field of quotients of
$\zz[\zz^n]$, $ F = \qq(t_1,\cdots,t_n)$; this was first observed by Le Dimet.  Here
we give a cohomological interpretation of this extension.  Our first application is
to prove that the representation is hermitian, extending a known result for braids. 
A simple proof of the concordance invariance of the represenation also follows.  The
cohomological approach leads to algorithms for computing the Gassner
representation, and these in turn yield connections between the Gassner
represention of a string link and the Alexander matrix of the link closure of that
string link; a new factorization result for Alexander polynomials follows.  A
random walk, or probablistic, approach to the Gassner representation is given,
extending previous work concerning the 1-variable Burau representation.  
We next show that by suitably normalizing, the Gassner matrix determines, and is
determined by, finite type link invariants. The paper
concludes with an interpretation of the determinant of the Gassner matrix in terms
of Reidemeister torsion, yielding an alternative approach to the factorization of
the Alexander polynomial.

\vfill
\eject

\begin{center}
{\bf \large THE GASSNER REPRESENTATION }
\end{center}

\begin{center}{\bf \large FOR STRING
LINKS}
\end{center}

\begin{center}{Paul Kirk, Charles Livingston, and Zhenghan Wang}
\end{center}
\begin{center}April 29, 1998\end{center}

\section{Introduction}

The correspondence between links and braids, as detailed by the Markov
theorem, has made  the study of  representations of the braid group a topic of
great interest in classical knot theory. Prior to the work of Jones
\cite{jones}, the most important such representation was the Burau
representation and its generalization, the Gassner representation, a
homomorphism defined on the  pure braid group on
$n$ strands,
$P_n$:
$$\gamma:P_n\to GL_n(\zz[\zz^n]).$$ A basic reference for the properties of 
the Burau and Gassner representations is \cite{birman}; more recent work
investigating the  Gassner representation includes
\cite{abdul,cochran, moody, penne}. 

In studying links, as opposed to knots, it has often proved useful to move
from the setting of braids to the more general objects of string links.  A
string link is simply a braid except that the strands are not required to be
monotonic; see Figure 1 and later discussions for details.  The usefulness of
this approach is perhaps best exemplified in the work of Habegger and Lin
\cite{habegger} in which a link homotopy classification is achieved through the
interplay between links and string links.  Part of the usefulness of braids
is that there is a natural group structure on the set of all braids, but there
is no similar group structure on the set of string links. However, once a
natural equivalence relation is placed on the set of string links, either
string link homotopy or concordance, a natural group structure reappears.

A significant advance in the study of links occurred in the work of Le Dimet
on link concordance.  In
\cite{ledimet} he observed that the Gassner representation extends to the
group  of concordance classes of  string links  provided one extends the range
to
$GL_n(F)$ where
$F$ denotes the quotient field of
$\zz[\zz^n]$.  It is the purpose of this article to explore this extension.  We
present simple new perspectives on its definition, its computation, and its
application to obtain both general results and interesting families of
examples.

A completely different approach to exploring the Gassner representation of
string links is being investigated by Ted Stanford.  His work is now in
preparation.

We thank Jim Davis , X. S. Lin, and Kent Orr for contributions to this work.  In
particular, Orr showed us the proof of Proposition 2.3.  Davis helped us develop
much of the work concerning torsion, presented in Section 10.

\vskip3ex

\noindent{\bf Contents}

In Section 2 we construct the Gassner matrix of a string link; the argument
uses a single straightforward cohomology computation.  The definition is
very simple and completely trivial in the case of braids. We were surprised
that this definition does not seem to appear explicitly in the literature,  even in
the context of braids or the Burau representation. 

The Burau and Gassner representations of braids are unitary with respect to a
particular Hermitian form.  This has been observed in \cite{abdul, squier, long}. 
In Section 3 we show that this symmetry extends to the case of the Gassner
representation for string links.  In addition, we show that the Hermitian form
corresponds to Poincar\`e duality on an appropriate cohomology group.

In Section 4 we use the well-known relationship between cohomology and the Fox
calculus to describe an algorithm for computing the Gassner representation of
a string link in terms of a presentation of the fundamental group of its
complement.  This also sets up some linear algebra which proves useful in
later sections.

In Section 5 we give the argument that the resulting extension of the Gassner
representation is a concordance and even I-equivalence invariant of the
string link.  I-equivalence is the equivalence relation defined by
non-locally flat concordance.  (It can be shown that two links or string links
are I-equivalent if and only if adding local knots to each can result in
concordant objects.)

Section 6 contains our main results, exploring a simple relationship between
the Gassner matrix of a string link and the Alexander matrix of its closure. 
We show that the Gassner matrix of a string link shares many formal properties
with the Alexander matrix of a link and exploit this to define the {\it
Alexander function of a string  link}.    In particular we prove that the
Alexander polynomial of the closure of a string link can be factored as a
product of a polynomial associated to the Gassner representation of the string
link and a second polynomial representing a Reidemeister torsion naturally
associated to the string link.  This torsion invariant vanishes for braids
and hence our result generalizes the classical result relating the Gassner
matrix of a braid and the Alexander polynomial of the associated link
\cite{birman}.  (See also \cite{penne}.) In Corollary \ref{irncymuo} we
show that the
 Alexander ideals $E_k(\hat{L})$  of the closure
$\hat{L}$ of a pure string link $L$   link vanish  if and only if the Gassner matrix
$\gamma(L)$ has the eigenvalue $1$ with multiplicity $k$.  In Theorem
\ref{alexp} we show that the  1-variable Alexander polynomials of
the link closure of a string link $L$,
$\hat{L}$, and the knot closure, $\widehat{LB}$,  are related by the
formula
$$ \bar{\Delta}_{\hat{L}}= \bar{\Delta}_{\widehat{LB}}\cdot
{{\det(I-\tilde{\beta}(L) )}\over { {\det(I-\tilde{\beta}(L)\tilde{\beta}(B) )}}}.$$
  The {\it knot closure,}  $\widehat{LB},$ is formed from the link closure by
banding the components together in a natural way, described explicitly in Section
6.  (The correction term
$
\det(I-\tilde{\beta}(L) )/\det(I-\tilde{\beta}(L)\tilde{\beta}(B) )$ is a string
link I-equivalence invariant of $L$.)

Section 7 presents examples.  For instance, we demonstrate that the subgroup
of the concordance group of string links formed by braids is not a normal
subgroup. 

In Section 8 we describe an alternative approach to the Gassner representation
via the notion of {\it random walks} on the string link.  This generalizes
similar work for the Burau representation developed in \cite{lin-wang}. 
Using this approach we are able to describe a construction  which, beginning
with almost any string link, produces an infinite family of string links, none
of which are concordant (or even I-equivalent) and all of which have the same
closure as the original string link. We also  indicate in this section 
that in some sense the Gassner representation is the most general
representation that can be attained via this random walk approach.

Section 9 discusses the connection between the Gassner representation and
finite type invariants.  In particular, by suitably normalizing it is seen
that the Gassner representation is made up of finite type invariants.  A
consequence is that it is determined by Milnor invariants of the string link,
although no explicit connection has been identified.

Finally, in Section 10 we interpret some of the previous work, especially the
product formula of Section 6, in terms of torsion invariants associated to a
string link and its associated closure.

\section{A cohomological definition of the Gassner representation}
 
We begin by introducing notation which will be in force throughout this
article.  Let $n$ be a positive integer;  usually $n$ denotes  the number of
components of a (string) link.  We will denote by
$\Lambda$  the integral group ring of
$\zz^n$ which will be thought of as the ring of Laurent polynomials in the
variables
$t_1,\cdots,t_n$, $\Lambda=\zz[t_1,t_1^{-1},\cdots, t_n,t_n^{-1}]$.  Its
fraction field, $\qq(t_t,\cdots ,t_n)$, is the field of rational
functions in the $t_i$ and will be denoted by
$F$.  We will occasionally use the localization $\Lambda_S$ of $\Lambda$
obtained by inverting all elements in the set, $S$, of Laurent polynomials
$f\in\Lambda$ satisfying
$f(1,\cdots, 1)=1$.    Thus
$\Lambda\subset
\Lambda_S\subset F$. Note that the units in $\Lambda$ are elements of the form
$\pm t_1^{a_1} \cdots t_n^{a_n}$.

The formal definition of  a string link is as follows.   Given a positive
integer
$n$, fix $n$ points in the interior of the 2-disk $p_1,\cdots, p_n$ , say
$p_i= (-1/(i+1), 0)$.   A {\it string link of n components} is a smooth, proper,
oriented 1-dimensional  submanifold of
$D^2\times [0,1]$ homeomorphic to the disjoint union of $n$ intervals such
that the initial point  of each interval  coincides with some
$p_{i}\times\{0\}$ and the endpoint coincides with $p_j\times\{1\}$. In the
figures the orientation is such that the components  run from the bottom to the top
of the diagrams. 

  Equivalently, one can think of a string link as a proper
 embedding of $n$ intervals, rather than as a submanifold.  Figure 1 illustrates a
3-stranded string link.


\vskip5ex
\begin{figure}
\centerline{\epsfxsize=1.3in\epsfbox{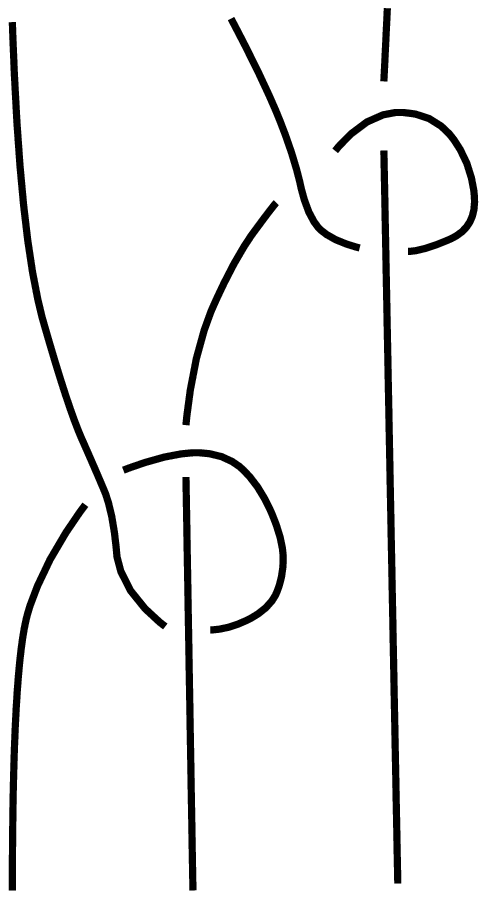}}
\caption{}
\end{figure}
\vskip5ex

 A {\it pure string link of n components} is a string link so that the
component labeled $i$ has initial point $p_i\times\{0\} $ and endpoint
$p_i\times\{1\}$.

Two string links are called {\it isotopic} if there is a smooth family of
string links interpolating between the two.

 Two string links are called
{\it concordant} if there is  a proper, locally flat, embedding  of
$n$ copies of
$[0,1]\times[0,1]$ in $D^2\times[0,1]\times [0,1]$   restricting to the
two string links when the last coordinate is $0$ or $1$, and mapping to
the proper $\{p_i\}\times\{ 0,1\}\times[0,1]$ when the middle coordinate is
$0$ or $1$.

Finally, two string links are called {\it I-equivalent} if we drop the
requirement that the concordance be locally flat.  

In this article we will simultaneously consider two situations which can be
roughly  formulated as the Gassner/multi-variable case and the
Burau/one-variable case.   The emphasis will be on the first case but
for the most part everything we say applies in either case.  However, the two
cases operate on slightly different categories of string links.  The proper
context for the Gassner/multi-variable case is to use {\it colored  string
links}; these are string links so that  each component is
indexed with an integer in $\{1,\cdots,n\}$.  For the Burau/single-variable
case the coloring is not needed or, equivalently, one can assume every
component is colored with the same color.

In particular, there is a multiplication on colored string links
$(L_1,L_2)\mapsto L_1\cdot L_2$ obtained by stacking $L_2$ on top of $L_1$.
This product  is defined only if the coloring  of the endpoints of $L_1$
coincide with the coloring of the initial points of $L_2$.  We will make some
attempt at stating things in this generality,  but will often restrict to
pure string links (with the component starting at
$p_i\times\{0\}$ colored with the label
$i$)  to simplify the exposition.

The {\it closure} of an $n$-component string link is the link of  
circles in $\rr^3$ obtained by joining the top points $p_i\times \{1\}$
with their corresponding bottom point $p_i\times \{0\}$ using $n$ arcs
in $\rr^3-(D^2\times I)$  which project to disjoint arcs in the plane.  For a
colored string link,  the closure is defined only if the colors ``match
up'', yielding a colored link of circles. See Figure 3.

We will now give  homological constructions of
the Burau and Gassner  representations.   These are surely known,
and of course follow easily from standard constructions of these
representations using the fact that the Fox calculus is merely a device for
encoding the differential on 1-cochains in the universal cover of a
space.  Our point of view is adapted from an article of Le Dimet's
\cite{ledimet} 

  Given a string link $L$ with $n$ components, let $X=D^2\times I -L$
denote the complement of the strands, and let $\pi=\pi_1 X$. The
abelianization of $\pi$ is isomorphic to $\zz^n$, and  the abelianization
map  
$$\epsilon:\pi\to
\zz^n=<t_1,\cdots,t_n>$$ is determined  by assigning to a  meridian  its
corresponding $t_i$. Multiplication by $\epsilon(\alpha)$ for
$\alpha\in\pi$ determines a local coefficient system on
$X$ with coefficients either $\Lambda$, $\Lambda_S$, or $F$, and hence
cohomology groups $H^*(X;\Lambda)$, $H^*(X;\Lambda_S)$, and  $H^*(X;F)$. 
Let $X_0= X\cap D^2\times\{0\}$ and $X_1= X\cap D^2\times\{1\}$; these
are both $n$-punctured discs and are canonically identified  (via the
homeomorphism $(x,0)\mapsto (x,1)$).  

Fix a point $p\in\partial D^2$ and let   $P\subset X$  be the arc  
$P=\{p\}\times I$ .
 
\begin{lemma}\label{asdfg}\par

\begin{enumerate}

\item $H^1(X_0;F)\cong F^{n-1}$ and $H^1(X_0,p;F)\cong F^n$.
\item The restriction maps $H^1(X;F)\to H^1(X_0;F)$ , $H^1(X;F)\to
H^1(X_1;F)$,\break $H^1(X,P;F)\to H^1(X_0,p;F)$, and $H^1(X,P;F)\to H^1(X_1,p;F)$
are all isomorphisms.
\end{enumerate}
\end{lemma}

We will prove Lemma \ref{asdfg}   shortly.   Assuming this, we can  
define the Gassner representation for a string link.
\begin{definition}

\begin{enumerate}
\item To a string link $L$ assign the automorphism
$$\gamma(L):H^1(X_0,p;F)\maptwo{\cong}{i_0^{-1}}H^1(X,P;F)
\maptwo{\cong}{i_1}H^1(X_1,p;F)\cong H^1(X_0,p;F).$$
(The last isomorphism being induced by the restriction of the  
map $D^2\times\{0\}\to D^2\times\{1\}$ given by $(z,0)\mapsto (z,1)$.)

 This
is called the {\it Gassner invariant} of
$L$.  It restricts to
  a  homomorphism on the semi-group of pure string links  (or the groupoid of
colored string links) 
$$\gamma:\call_n\to GL(H^1(X_0,p;F))\cong GL_n(F)$$ called the {\it Gassner
representation}.

\item To a string link $L$ assign the automorphism
$$H^1(X_0;F)\maptwo{\cong}{i_0^{-1}}H^1(X;F)
\maptwo{\cong}{i_1}H^1(X_1;F)\cong H^1(X_0;F).$$
This is called the {\it reduced Gassner invariant} and restricts to 
a  homomorphism on the semi-group of pure string links 
$$\tilde{\gamma}:\call_n\to GL(H^1(X_0;F))\cong GL_{n-1}(F)$$
called the {\it reduced Gassner representation}.\end{enumerate}
In these definitions the maps $i_0$ and $i_1$ are induced by the
inclusions $X_0\subset X$ and $X_1\subset X$.
\end{definition}

From this construction it is obvious that if $L_1$ and $L_2$ are two
colored string links and $L_1\cdot L_2 $ denotes the product of the two,
then 
$$\gamma(L_1\cdot L_2)=\gamma(L_1)\gamma(L_2).$$

Similarly define the {\it Burau representation},
$$\beta(L) \in GL_n(\qq(t)),$$ by composing $\gamma$ with the ring homomorphism
$F\to \qq(t)$ taking each $t_i$ to $t$.  The Burau
representation is multiplicative  on all string links, not just pure
string links.  We will show below that these representations are string
link concordance invariants, and  thus descend to homomorphisms on the
concordance groups (this was first observed in \cite{ledimet}.)

   Moreover, there is a canonical choice of basis of
$H^1(X_0,p;F)$ so that one can identify these with matrix
representations.  It will
follow as an elementary exercise that  if the string link is a (pure)
braid, then the resulting matrix representation is just the standard
Gassner representation.

We have used cohomology to define the invariants $\gamma$ and $\beta$.  We could just as well have
used homology, and Theorem \ref{hgecg} shows that the ``homology Gassner'' and ``cohomology
Gassner'' are equivalent representations.

The proof of the following result was
pointed out to us by K. Orr.

\begin{proposition}\label{ssalkdja}
Let $(X,Y)$ be a pair of path connected CW complexes and $\epsilon:\pi_1X\to
\zz^n$ a homomorphism.  Consider the corresponding $F=\qq(t_1,\cdots,t_n)$
coefficient system on the pair $(X,Y)$.

 Suppose that the inclusion of
$Y$ in
$X$ induces an isomorphism on homology with (untwisted) $\qq$
coefficients.  Then $H^*(X,Y;F)=0$. (In fact $H^*(X,Y;\Lambda_S)=0$.)
\end{proposition}
\noindent{\sl Proof.}
Let  $(C_*(X,Y),\partial)$ denote the cellular chain complex with
$\qq$ coefficients for the pair $(X,Y)$ and let
$(C_*(\tilde{X},\tilde{Y}),\tilde\partial)$ denote the chain complex of
the covering of the  pair determined by $\epsilon$.  Fix lifts of the cells of
$(X,Y)$  to $(\tilde{X},\tilde{Y})$ to  get a free $\Lambda$ basis of
$C_*(\tilde{X},\tilde{Y})$.

By hypothesis, $C_*(X,Y)$
is acyclic.  Thus there exists a chain contraction $D:C_*(X,Y)\to
C_*(X,Y)$, i.e.   a degree 1 map  satisfying $D\partial+\partial
D=Id$.   Using the chosen free $\Lambda$ basis for
$C_*(\tilde{X},\tilde{Y})$ and the formula for $D$ one can define a
degree 1 map (in fact a chain homotopy) $\tilde{D}:
C_*(\tilde{X},\tilde{Y})\to
 C_*(\tilde{X},\tilde{Y})$. Explicitly,  if $De=\sum_i q_i f_i$ define
$\tilde{D}\tilde{e}=\sum_i q_i\tilde{f}_i$ where
$\tilde{e},\tilde{f}_i$ are the chosen lifts of $e, f_i$.

By construction
$\Phi=\tilde{D}\tilde{\partial}+\tilde{\partial}\tilde{D}$ is a chain
map whose matrix in the chosen basis augments to the identity map.
That is,  if $a:\Lambda\to \qq$ is the augmentation $t_i\mapsto 1$, then
$a(\Phi)=Id$.    Dualizing, the induced chain homotopy on the cochain
complex
Hom$_{\Lambda}(C_*(\tilde{X},\tilde{Y}),F)$ is a
chain homotopy from $\Phi^*$ to $0$; hence $\Phi^*$ induces the zero map on
the cohomology $H^*(\tilde{X},\tilde{Y};F)$.  But
since the determinant of $\Phi^*$  is a non-zero  element of
$\Lambda\subset
F$ (it augments to $1$), $\Phi^*$ is a chain
isomorphism. In other words, the zero map on cohomology is an
isomorphism, and thus the cohomology must vanish.  (Notice that the argument
works with $\Lambda_S$ coefficients, since the determinant of $\Phi^*$
augments to $1$, and hence this  determinant is in $S$.)
\qed

Lemma \ref{asdfg} follows from Proposition \ref{ssalkdja} by observing
that the inclusions $X_0\subset X$ and $X_1\subset X$ induce
isomorphisms in (untwisted) rational homology, using for example
Alexander duality, and carrying out the (simple) calculations for $X_0$
or $X_1$ which are homotopy equivalent to a wedge of circles (this is
also done below).  

 We summarize these observations in the following theorem.

\begin{theorem}
The assignment $L\mapsto \gamma(L)$ defines an isotopy
invariant of the $n$-component string link $L$ with values in
$GL_n(F)$.  It is multiplicative under the  
multiplication of colored string links obtained by stacking one above the other:
$$\gamma(L_1\cdot L_2)=\gamma(L_1)\cdot
\gamma(L_2).$$
Restricting to    pure braids   yields the classical Gassner
representation.
\end{theorem}

\section{Symmetry of the Gassner representation}
It has been observed by several authors \cite{abdul, squier, long} that the 
reduced Burau and Gassner representations of braids are unitary with respect to an
appropriate Hermitian form.  Here we will present a proof of this fact that
generalizes to the Gassner representation for string links.  However, rather than
work in specific coordinates we will see that the symmetry is a consequence of
Poincar\'e duality, and that the representation is unitary with respect to a skew
symmetric form on
$H^1(X;F)$ corresponding to the cup product.  

Before presenting the argument, we  briefly review  Poincar\'e duality and
intersection forms in the setting of local coefficients.  To simplify
notation we will describe the case of closed manifolds.  The extension to
bounded manifolds (which we will use) is straightforward.  We will also
restrict to coefficients in a field, $F$, with involution, acted on by
$\pi$ via an involution preserving representation $\rho: Z[\pi]
\rightarrow F$.

Let $X$ be a triangulation of a closed oriented $n$-manifold $M$ and let
$X'$ denote the dual cell decomposition.  We again denote the fundamental
group by $\pi$ and assume $\pi$ acts on the right on the universal cover
$\tilde{M}$ and on the left on $F$ via $\rho$.  We denote by
$\tilde{X}$ and $\tilde{X}'$ the corresponding decompositions of
$\tilde{M}$.  Hence, the chain complexes $C_*(\tilde{X}) \otimes_{\pi} F$ and
 Hom$_{\pi}(C_*(\tilde{X}),F)$ determine the twisted homology and cohomology
groups.  Note that  Hom$_{\pi}$ denotes the set of left $\pi$ homomorphisms;
$\pi$ acts on $ C_*(\tilde{X})$ on the left via right multiplication by
inverses.  Both $H_*(M,F)$ and $H^*(M,F)$ are right
$F$ vector spaces.

There is an intersection form $\langle\ ,\ \rangle :  C_i(\tilde{X}) \times
C_{n-i}(\tilde{X'}) \rightarrow Z[\pi]$ defined by $\langle a , b\rangle  =
\sum_g{(a \cdot bg)g}$.  This in turn determines a pairing
 $$\langle\ ,\ \rangle :  H_i(M;F) \times H_{n-i}(M;F) \rightarrow F$$
satisfying 
\begin{equation}\label{unitary}
\langle af,b\rangle  = \langle a,b\rangle f,\ \ \ \langle
a,bf \rangle  = f^{-1} \langle a,b \rangle \ \ \  and \ \
\ \langle a,b\rangle  = -\overline{\langle b,a\rangle }.\end{equation} 
The map
$\Phi:C_i(\tilde{X})
\rightarrow \mbox{Hom}(C_{n-i}(\tilde{X'}),F)$ is defined by $\Phi(a)(b) = \langle
a, b \rangle$; $\Phi$  induces the Poincar\'e isomorphism, $\Phi:H_i(M,F)
\mapright{\cong} H^{n-i}(M;F)$.  Note too that $\Phi$ induces an isomorphism
 Hom$_F(H_i(M;F),F)$, where  Hom$_F$  denotes left $F$ homomorphisms and
$H_i(M;F)$ is a left $F$ vector space via right multiplication by inverses.

To apply this to the Gassner representation, consider the diagram below. 

\[
\begin{diagram}\dgARROWLENGTH=2ex
\node{H^1(X_0;F)}
    \node{H^1(X;F)}\arrow{w}\arrow{e}
    \node{H^1(X_1;F)}\\
\node[2]{H_2(X,X_0 \cup X_1;F)}\arrow{sw}\arrow{se}\arrow{n,l}{\Phi}\\ 
    \node{H_1(X_0;F)}\arrow{e}\arrow[2]{n,l}{\Phi} 
    \node{H_1(X;F)}
 \node{H_1(X_1;F)}\arrow{w}\arrow[2]{n,l}{\Phi}
\end{diagram}\]

Every map in this diagram is an isomorphism. The
top row yields a homomorphism $\nu_c:H^1(X_0,F)
\rightarrow H^1(X_1;F)$ and the bottom row yields  $\nu_h:H_1(X_0;F)
\rightarrow H_1(X_1;F)$. Composing either with the map induced by identifying
$X_0$ with $X_1$ yields the reduced Gassner and {\it homology} Gassner
representations.  Commutativity of the upper two squares follows from the
naturality of Poincar\'e duality.  That the bottom central square commutes
follows readily from the associated long exact sequence.  Hence:

\begin{theorem}\label{hgecg} $\Phi \nu_h = \nu_c \Phi$.  In particular the homology
and cohomology Gassner representations are equivalent. \qed \end{theorem}

That the homology Gassner
representation is unitary with respect to the intersection form follows from the next
result.

\begin{theorem} For $a  \in H_1(X_0;F)$ and $b \in H_1(X_1;F)$,
$\langle \nu_h(a), b\rangle  = \langle a, \nu_h^{-1}(b)\rangle $.\end{theorem} 

\noindent{\sl Proof.}  Given the definition of $\nu_h$ it is clear that there is a class
$A \in H_2(X,X_0 \cup X_1,F) $ with $\partial(A) = a -
\nu_h(a)$.  Hence, $\langle \nu_h(a),b\rangle  = \langle A,b\rangle  =
\langle A,\nu_h^{-1}(b)\rangle  =
\langle a,\nu_h^{-1}(b)\rangle  $. (Here the middle two intersections are on
the relevant relative homology groups.)\qed

\begin{corollary}\label{oremvcxmgigt}
The reduced homology  Gassner invariant of a string  link $\tilde{\gamma}(L)$
satisfies
$$\langle  \tilde{\gamma}(L)(a),\tilde{\gamma}(L)(b)\rangle=\langle
a,b\rangle.$$
\end{corollary}
\qed

In particular one obtains the following generalization to string links of the main
result of
\cite{abdul}.   

\begin{theorem}\label{traceresult}
There exists an imbedding $e:F\hookrightarrow \cc$ of fields so that 
$e\circ \tilde{\gamma}(L)$ is a unitary complex matrix for all string links $L$. In particular the 
trace Tr$({\gamma}(L))$ equals $n$ if and only if $\gamma(L)$ is the identity matrix.  The same
statements apply to the Burau representation.
\end{theorem}
\noindent{\sl Proof.}  Given  unit complex numbers $z_1,\cdots,z_n$ so that the
transcendence degree of $\qq[z_1,\cdots,z_n]$ over $\qq$ is $n$, one obtains a conjugation
preserving imbedding
$j:F\to
\cc$   fields by taking $t_i$ to $z_i$.     A straightforward geometric calculation shows
that if  one chooses $z_j= e^{2\pi i a_j}$ with $0<a_j<{1\over{n+1}}$, then the resulting
Hermitian (over $\cc$) intersection pairing
\begin{equation}\label{rotori}
H^1(X_0;\cc_\alpha)\times H^1(X_0;\cc_\alpha)\to \cc
\end{equation}
is  definite.  (More precisely, after multiplying by $i$ the skew-hermitian pairing 
(\ref{rotori}) becomes   definite.)   Here
$\alpha$ denotes the induced
$\cc$ local coefficients. 

 (A
``geometry-free'' argument can be obtained using the results of
\cite{aps2}.  The signature of this pairing is computed by formula
$$\mbox{Sign}(X_0)-\mbox{Sign}_\alpha(X_0)=\rho_{\alpha}(\partial
X_0)=-(1-2\sum_{i=1}^n a_i) +
\sum_{i=1}^{n} (1-2a_i)= n-1$$
provided that $0<\sum_{i=1}^n a_i<1$.  Since
$\mbox{Sign}(X_0)=0$ and $H^1(X_0;\cc_\alpha)=\cc^{n-1}$, the  pairing
(\ref{rotori}) is negative definite after multiplying by $i$.)

 Fix a basis for $H^1(X_0;F)$.  In this basis the intersection   $\langle a, b \rangle$
equals
$b^*Ja$ for some matrix $J$ which satisfies $J^*=-J$ by Equation \ref{unitary}.   Thus the
complex matrix $K=i\cdot j(J)$ is definite and satisfies $K^*=K$. Therefore there exists
an isomorphism
$\psi:\cc^n\to \cc^n$ taking $K$ to $ \pm I$.    Let $e(f)=\psi^{-1} j(f) \psi$.  

In the chosen basis, Corollary \ref{oremvcxmgigt} says that   for any string
link $L$
$$\tilde{\gamma}(L)^*J\tilde{\gamma}(L)=J.$$   Applying $e$ one sees that
$$e(\tilde{\gamma}(L))^* \ iI \ e(\tilde{\gamma}(L))= iI$$
and since $iI$ is central 
$$e(\tilde{\gamma}(L))^* \  e(\tilde{\gamma}(L))= I,$$
so that $e(\tilde{\gamma}(L))$ is unitary, as claimed.     

Since $e$ is an imbedding, the trace of  $\tilde{\gamma}(L)$ is $n-1$ if and only if the
trace of
$e(\tilde{\gamma}(L))$ is $n-1$,  which occurs only if $e(\tilde{\gamma}(L))=I$ since this
matrix is unitary.  We will show below (Proposition \ref{sadpoi}) that $\tilde{\gamma}$ is
the reduction of
$\gamma$ with respect to a  fixed $1$-eigenspace for $\gamma$, so that the trace of
$\gamma$ equals Tr$(\tilde{\gamma})+1$.  Thus the trace of $\gamma(L)$ equals $n$ if and
only if the trace of
$\tilde{\gamma}(L)$ equals $n-1$.

All the arguments in this section carry out without significant change to the 
1-variable Burau case.\qed

\section{The Gassner
representation, group cohomology, and the Fox calculus}\label{wekvoi}

We will next 
give a different argument for Lemma \ref{asdfg} using the Fox calculus.
This will make more transparent how easy this invariant is to calculate
and will introduce the Wirtinger--Fox matrix of a string link.  It will be
seen that this matrix carries more information about the link than the
cohomology alone.   Moreover,  we will describe a canonical basis of
$H^1(X,P;F)$ which we use to identify $\gamma$ with a
matrix representation. This can be used to 
make explicit the relation of the Gassner representation described
above to the classical Gassner representation, and to the results of
\cite{lin-wang}.  Finally, in twisted generalizations of the Gassner
representation that we will introduce in a later article,  the hypotheses of
Proposition
\ref{ssalkdja} need not hold and it  will be convenient to have an explicit
model to compute with.

The model  we will use to compute cohomology
is the  reduced  bar resolution of the fundamental
group. Since the cohomology of a space and its
fundamental group coincide in degrees
$0$ and $1$ (and since in any case $X, X_1$ and $X_0$ are aspherical)
and these are the only groups entering into the definition of the
Gassner representation, we lose nothing  by computing in group
cohomology.

 Recall that to a group
$\pi$ and a left
$\zz[\pi]$ module $M$  the reduced bar resolution assigns the cochain
complex
$(C^*(\pi,M),d)$ where
$$C^n(\pi,M)=\mbox{Maps}(\underbrace{\pi\times\cdots\times\pi}_{n\
\hbox{times}},M)$$
and the differentials are given by standard formulas derived from the
boundary operators on simplices.  Germane for us will be the formulas
$$d^1:C^1(\pi,M)\to C^2(\pi,M), \ \ \  d^1(f)(x,y)= f(x) + x\cdot f(y) -
f(xy)$$
and
$$d^0:C^0(\pi,m)=M\to C^1(\pi,M), \ \ \ d^0(m)(x)= m- x\cdot m.$$

Given a string link, let $\pi$ denote the fundamental group of its
complement.    A choice of projection of the string link
determines a Wirtinger presentation of its complement. Fix a projection
of the string link of $n$ components and assume it has $c$ crossings. It
turns out to be convenient (though not essential) to assume that the
projection has the property that every component of the string link projects
to a curve   which  passes   under
at some crossing.   This is automatically satisfied  if the
closure is not a split link. Notice that this implies that $c\ge n$.  
One quick way to ensure this is to add a small ``kink'' near the top endpoints
of the strands  (see Figure 2).

We will label the Wirtinger generators
in the following way.  Let
$\mu_1,\cdots,
\mu_n$ denote the generators corresponding to meridians which lie in
$D^2\times \{0\}$, $\mu_1',\cdots, \mu_n'$ the generators corresponding
to the meridians which lie in $D^2\times \{1\}$, and $z_1\cdots z_{c-n}$
denote the remaining generators.  Figure  2 illustrates the labeling
scheme.   Moreover we will always order this basis as 
$$\mu_1, \cdots,\mu_n, z_1, \cdots, z_{c-n}, \mu_1' , \cdots , \mu_n'.$$
 Each relation in the
Wirtinger presentation is of the form
$\alpha\beta\alpha^{-1}\gamma^{-1}$ where
$\alpha,\beta$ and
$\gamma\in\{\mu_i,z_i,\mu_i'\}$.  
For example, in the following figure, the middle right crossing
determines the Wirtinger relation  $\mu_2 z_3\mu_2^{-1}z_2^{-1}$  and
the top left crossing  determines the Wirtinger relation  $z_1 \mu_1'
z_1^{-1}z_1^{-1}$ (which of course reduces to $\mu_1'z_1^{-1}$).


\vskip5ex
\begin{figure}
\centerline{\epsfxsize=1.5in\epsfbox{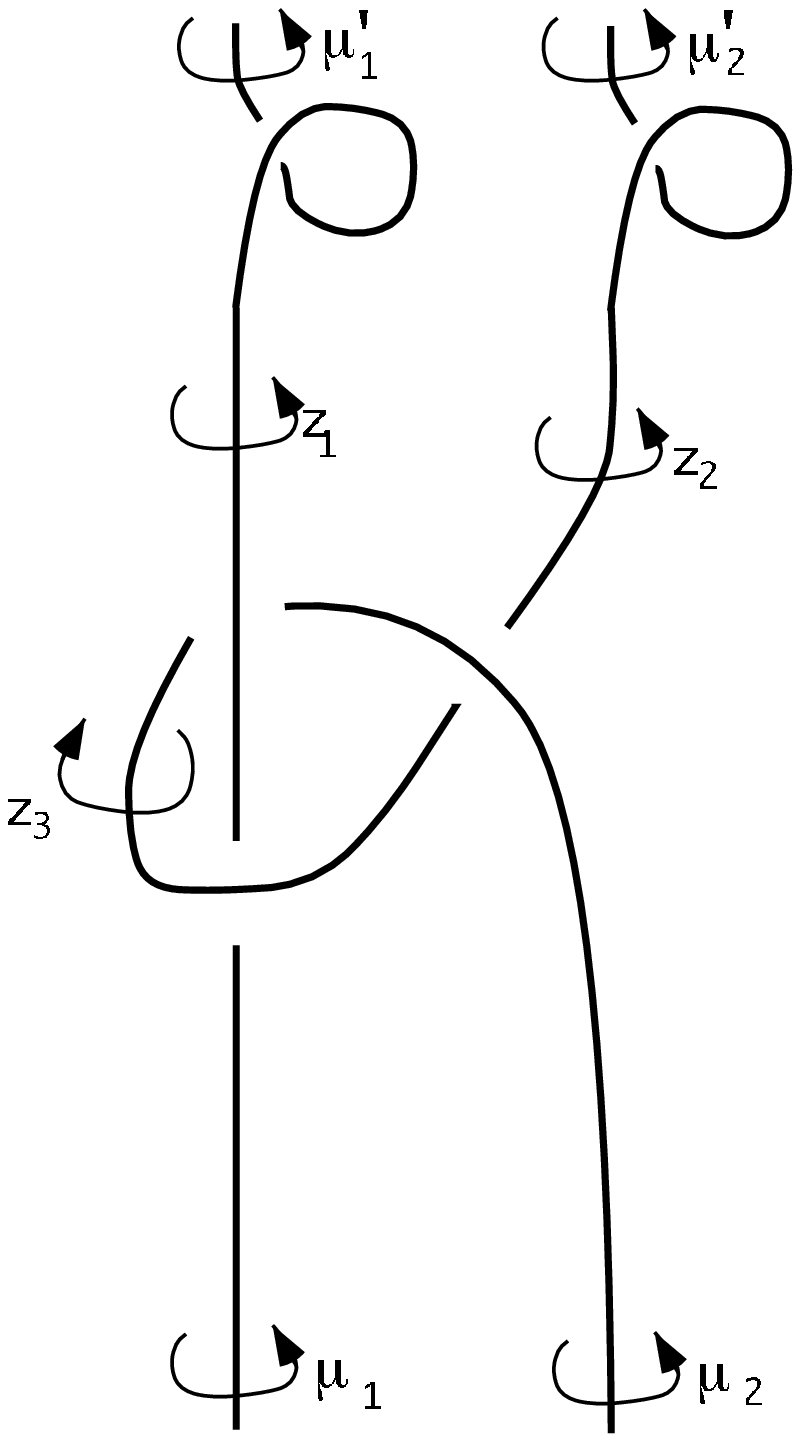}}
\caption{}
\end{figure}
\vskip5ex

We will use the Wirtinger presentation of the fundamental group of
the string link to compute   cohomology, using the following lemma.
This simple  lemma 
forms the basis for the Fox calculus and its relationship to group cohomology. 
Any statement made in this article about the Fox calculus is derived from this
lemma together with the definition of a 1-cocycle.

\begin{lemma}\label{xovci}
Suppose that a presentation $ <x_i, i\in I | r_j, j\in J>$ is given
for the group $\pi$, and an action $\pi\to \mbox{Aut}(M)$ is given.  Then 
\begin{enumerate}
\item A 1-cocycle with values in $M$ is determined by its values on the
$x_i$. 

\item Given a collection of elements $m_i\in M, i\in I$, there exists a
(unique) 1-cocycle $f$ on the free group generated by  the $x_i$
satisfying $f(x_i)=m_i$ (with respect to the composite action $<x_i>\to
\pi\to\mbox{Aut}(M)$).

\item The cocycle $f$ given in the previous assertion well-defines a
1-cocycle on $\pi$ if and only if $f(r_j)=0$ for all $j$.
\end{enumerate}
\end{lemma} \qed

Fix a $\pi$-module $M$ determined by an action  $e:\pi\to
\mbox{Aut}(M)$.

  Denote by $K$ the free
group on $\{\mu_i,z_i,\mu_i'\}$.  Let
$f\in Z^1(K;M)$ be some 1-cocycle on $K$, which, according to Lemma
\ref{xovci}, is uniquely determined by the values of $f$ on the free
generators. Then using the fact that
$f$ is a cocycle, i.e. a crossed homomorphism on the free group,  one
computes:
\begin{eqnarray}\label{lksdj}
 f(\alpha\beta\alpha^{-1}\gamma^{-1})&=& f(\alpha)+e(\alpha)
f(\beta) - e(\beta)f(\alpha) - f(\gamma)\\
&=& (1-e(\beta))f(\alpha) + e(\alpha)f(\beta) -
f(\gamma)\nonumber
\end{eqnarray}

  Lemma \ref{xovci} implies that $f$ determines a 1-cocycle on 
$\pi$ if and only if for each
Wirtinger relation $\alpha\beta\alpha^{-1}\gamma^{-1}$ one has
 $f(\alpha\beta\alpha^{-1}\gamma^{-1})=0$.

 We use the
Fox calculus to reduce the computation of cocycles to linear algebra.
Restrict to
$M=F$ and
$e=\epsilon$. Notice that in Equation \ref{lksdj} the elements
$e(\alpha), e(\beta)$, and $e(\gamma)$ are in the set $\{t_i\}_{i=1}^n
$.
  Thus each Wirtinger relation determines a row in a
$c\times (n+c)$   matrix  using Equation \ref{lksdj} in the following way:
\begin{enumerate}
\item The columns
correspond to the (ordered) $(n+c)$-tuple
$$(\mu_1,\cdots,\mu_n,z_1,\cdots,z_{c-n},\mu_1'\cdots\mu_n').$$   
\item     The
rows correspond to the relations, so the row corresponding to the relation
$\alpha\beta\alpha^{-1}\gamma^{-1}$ has
$1-e(\beta)$ in the column corresponding to $\alpha$, $e(\alpha)$ in
the column corresponding to $\beta$, and $-1$ in the column
corresponding to $\gamma$.  
\end{enumerate}

 For
example,   a crossing with Wirtinger  relation $\mu_i z_j\mu_i^{-1}
(\mu_k')^{-1}$   contributes a row  with $1-t_j$ in the $i$th column,
$t_i$ in the $(n+j)$th column, and
$-1$ in the $(c+k)$th column.  Write this matrix $\pmatrix{A&B&C\cr}$
where $A$ is a $c\times n$ matrix, $B$ is $c\times (c-n)$, and $C$ is
$c\times n$.

From Lemma \ref{xovci} one immediately concludes the following.
 \begin{lemma}\label{paslapsl}
Given a collection 
$\lambda_1,\cdots, \lambda_n,$  $\tau_1, \cdots, \tau_{c-n},
\lambda_1',\cdots
\lambda_n'\in F$,  there is a cocycle $f\in Z^1(\pi;F)$ satisfying 
$f(\mu_i)= \lambda_i$, $f(z_i)= \tau_i$, and  $f(\mu_i')= \lambda_i'$ if
and only if 
\begin{equation}\label{alexander}
\pmatrix{A&B&C}
\pmatrix{\lambda_1\cr\vdots\cr \lambda_n\cr \tau_1\cr \vdots\cr
\tau_{c-n}\cr
\lambda_1'\cr\vdots\cr
\lambda_n'\cr}=0\end{equation}
or, equivalently
\begin{equation}\label{sdlkas} 
\pmatrix{A&B}
\pmatrix{\lambda_1\cr\vdots\cr \lambda_n\cr \tau_1\cr \vdots\cr
\tau_{c-n}\cr}=-C \pmatrix{
\lambda_1'\cr\vdots\cr
\lambda_n'\cr}\end{equation}
\end{lemma}
\qed
   The matrix $\pmatrix{A&B&C\cr}$ will be important in what
follows.  We will call it the {\it Wirtinger--Fox matrix of the string link}.
The next lemma substitutes for Proposition \ref{ssalkdja} in the present context.

\begin{lemma}\label{paorjvc}
The square  matrix $\pmatrix{A&B\cr}$ has entries in
$\Lambda=\zz[t_1,t_1^{-1},\cdots ,t_n,t_n^{-1}]$.  Its determinant is a Laurent
polynomial in the
$t_i$ whose value at $1$ equals
$\pm1$.   Therefore
$\pmatrix{A&B\cr}$ is invertible over $\Lambda_S$ and hence also over $F$. 
\end{lemma}
\noindent{\sl Proof.} From its definition $\pmatrix{A&B\cr}$ has entries in
$\Lambda$.  If we set all the
$t_i$ equal to $1$, the equation 
$$  (1-\epsilon(\beta))f(\alpha) + \epsilon(\alpha)f(\beta) -
f(\gamma)=0$$ reduces
to   \begin{equation}\label{xzpob} 
f(\gamma)=f(\beta) \end{equation}
at each crossing.   This implies that if
$a:\Lambda\to
\zz$ denotes the augmentation
$t_i\mapsto 1$, then the system of linear equations over $\zz$
$$a \pmatrix{A&B} 
\pmatrix{x_1\cr\vdots\cr x_n\cr y_1\cr \vdots\cr
y_{c-n}\cr}=-a(C) \pmatrix{
s_1\cr\vdots\cr
s_n\cr}$$
has a unique integer solution for any choice of $s_1,\cdots , s_n$. 
One way to see this explicitly is to number the crossings as follows: let the
first crossing  be the  crossing  at which the
strand with Wirtinger generator 
$\mu_1$ ends, the second crossing is the crossing at which the strand with
Wirtinger generator $\mu_2$ ends, and so forth, going through the $\mu_i$ and
then the $z_i$ in order.  With this choice of ordering of the generators, the
matrix $a(A\ B)$ is upper triangular with 1s on the diagonal. 
 Therefore the determinant of 
$a
(A\ B)
$  equals
$1$.  The lemma follows.
\qed

  To
simplify and clarify notation, we define, for   $Z=X, X_1, X_0, P$, or $p$ 
$$C^i(Z)=C^i(\pi_1(Z);F).$$
 Since space and group cohomology agree for connected spaces in dimensions $0$ and
$1$,
$H^i(Z;F)$ is the cohomology of $C^*(Z)$ for these spaces for $i=0$ and
$1$ (and even for all $i$ since each $Z$ is aspherical).  

Also, define $C^i(X,P)=\mbox{Ker} \ C^i(X)\to C^i(P)$ and 
$C^i(X_1,p)=\mbox{Ker}\ C^i(X_1)\to C^i(p)$.  
We use similar notation for cocycles ($Z^i(Z)$) and coboundaries
($B^i(Z)$).

The  inclusion  $i_1: (X_1,p)\to (X,P)$ induces maps on cochains.  Since
$\pi_1X_1$ is free on the $\mu_i'$, Lemma \ref{xovci} implies that
$Z^1(X_1)\cong F^n$ and that an identification  is given by the assignment
$f\mapsto (f(\mu_1'),\cdots, f(\mu_n'))$.   On the other hand, Lemma 
\ref{paorjvc} and Equation \ref{sdlkas} show that the restriction map 
on 1-cocycles
$$i_1:Z^1(X)\to Z^1(X_1)$$
taking $f$ to its restriction to $\pi_1(X_1)$ is an isomorphism.  

Note that $d^0:F=C^0(X)\to C^1(X)$ is given by $d^0(v)(x)=
(1-\epsilon(x))v$.  Since, for example, $1-\epsilon(\mu_i')= 1-t_i\ne 0$, $d^0$
is non-zero and hence injective.  Similarly  $d^0:F=C^0(X_1)\to C^1(X_1)$
is injective and so by naturality we conclude that
$$i_1:H^1(X;F)\to H^1(X_1;F)$$
is an isomorphism of $F$-vector spaces of dimension $n-1$.
A similar argument applies to show that $i_0:H^1(X;F)\to H^1(X_0;F)$ 
is an isomorphism.
 
We claim that every map in the diagram
\[
\begin{diagram}
\node{H^1(X,P;F)}\arrow{s,l}{i_1}
    \node{Z^1(X,P)}\arrow{w}\arrow{e}\arrow{s,l}{i_1}
    \node{Z^1(X)}\arrow{s,l}{i_1}\\
\node{H^1(X_1,p;F)} 
    \node{Z^1(X_1,p)}\arrow{w}\arrow{e} 
    \node{Z^1(X_1)}
\end{diagram}\]
is an isomorphism.

Indeed, since $C^0(X,P)= \mbox{Ker} \ C^0(X)\to C^0(p)= \mbox{Ker} \
Id:F\to F$, it follows that the natural projection  $Z^1(X,P)\to
H^1(X,P;F)$ is an isomorphism. Similarly
$Z^1(X_1,p)\to H^1(X_1,p;F)$ is an isomorphism.   Since $P$ is an arc,
$\pi_1(P)=(1)$, so $C^i(P)\cong F$ for all $i$.  The restriction map
$C^1(X)\to  C^1(P)$ takes a 1-cochain to its value at $1\in \pi_1X$. 
Since all 1-cocycles vanish on the identity element,  
$Z^1(X)\subset C^1(X,P)$.  But since $Z^1(X,P)= Z^1(X)\cap C^1(X,P)$  we
  conclude that $F^n=Z^1(X)= Z^1(X,P)= H^1(X,P;F)$.  A similar argument
applies to $(X_1,p)$.  Finally the five lemma shows that the left
vertical map in the diagram is an isomorphism, and so all maps in the diagram are
isomorphisms.

This shows how to interpret $\gamma(L)$ as a matrix.  Just 
define
$\gamma(L)$  to be the matrix representing the composition of isomorphisms
\[ 
\begin{diagram}\dgARROWLENGTH=2.8ex
\node{F^n}\node{H^1(X_0,p;F)}\arrow{w,t}{\alpha_0}
\node{H^1(X,P;F)}\arrow{w,t}{i_0}\arrow{e,t}{i_1}\node{H^1(X_1,p;F)}
\arrow{e,t}{\alpha_1}\node{F^n}
\end{diagram}
\]
The isomorphism $\alpha_0$ takes a cocycle $f$ to the vector
$(f(\mu_1)\cdots,f(\mu_n))$ and the isomorphism $\alpha_1$ takes a cocycle
$f$ to the vector
$(f(\mu_1')\cdots,f(\mu_n'))$.
 
\vskip2ex

The following observation will be important in a later section.   Using
Equation \ref{sdlkas} and the definition of $\gamma(L)$, for each
$i=1,\cdots,n$  there is a (unique)  vector $(z_{1i},\cdots,z_{c-n\ i})$
so that 
$$\pmatrix{A&B\cr}
\pmatrix{\gamma(L)_{1i}\cr\vdots\cr\gamma(L)_{ni}\cr
z_{1i}\cr\vdots\cr z_{c-n\ i}\cr}=
-C\pmatrix{0\cr\vdots\cr1\cr\vdots\cr0\cr}.            
$$

This implies the following.

\begin{proposition}\label{xys} The Wirtinger--Fox matrix $(A\  B \ C)$ of a
projection of a string link $L$ and the Gassner matrix $\gamma(L)$ are related by
the formula 
\begin{equation} \label{sdajasd}
\pmatrix{A&B\cr}\pmatrix{\gamma(L)\cr Z\cr}=-C
\end{equation}
where $Z$ is some matrix that depends on the projection. \qed
\end{proposition}

To see that $\gamma(L)$ gives the classical Gassner representation when
restricted to pure braids, one can either check by hand that this
construction gives the correct value on the standard generators of the
pure braid group, or else one can observe that the standard definition
in terms of  free group automorphisms gives rise via cohomology (or
equivalently the Fox calculus)   to the same $\gamma(L)$.  We leave
the details to the reader.

\section{Concordance and boundary links}

In this section we prove that $\gamma(L)$ is an I-equivalence
invariant of $L$ and vanishes on boundary links.  Let $\calc_n$ denote the
group of I-equivalence classes of
pure $n$-component string links. 

\begin{theorem}\label{iequivalence}
 The assignment to a string link $L$ the matrix $\gamma(L)$ defines an
I-equivalence invariant of $L$.  Moreover the induced function   
  on the group of I-equivalence classes
of
pure $n$-component string links 
$$ {\gamma}:\calc_n\to GL_n(F)
$$
is a group homomorphism.
\end{theorem}
\noindent{\sl Proof.}
Suppose that $L$ and $L'$ are I-equivalent string links.  Let $X$ and
$X'$ denote the complements of $L$ and $L'$ in $D^2\times I$, and let
$Y$ denote the complement of an embedding of $n$ copies of $I\times I$
in $D^2\times I \times I$ exhibiting the concordance.  Thus, using the
obvious notation, we have a commutative diagram of inclusions:
\[
 \begin{diagram}
\node{X_0}\arrow{e}\arrow{s}\node{Y_0}\arrow{s}\node{X_0'}\arrow{w}\arrow{s}\\
\node{X} \arrow{e}\node{Y}\node{X'}\arrow{w}\\
\node{X_0}\arrow{n}\arrow{e}\node{Y_0}\arrow{n}\node{X_0'}\arrow{w}\arrow{n}
\end{diagram}
\]

Applying cohomology with twisted coefficients in $F=\qq(t_1,\cdots,
t_n)$ to this diagram and invoking Proposition
\ref{ssalkdja}
we see that in the following commutative diagram all maps are
isomorphisms. 
\[
 \begin{diagram}\dgARROWLENGTH=2.3ex
\node{H^1(X_0,p)}\node{H^1(Y_0,p\times
I)}\arrow{w}\arrow{e}\node{H^1(X_0',p)}\\
\node{H^1(X,P)}\arrow{n} \arrow{s} \node{H^1(Y,P\times
I)}\arrow{n}\arrow{s}\arrow{w}\arrow{e}\node{H^1(X',P)}\arrow{n}\arrow{s}\\
\node{H^1(X_0,p)} \node{H^1(Y_0,p\times I)}\arrow{e}\arrow{w}
\node{H^1(X_0',p)} 
\end{diagram}
\]
Every cohomology group in the top and bottom rows are canonically
identified with
$F^n$ as described in the previous section.  With these identifications   the two top
horizontal and bottom two horizontal isomorphisms are given by the
identity map.    Since
$\gamma(L)$ is given by the composite of the two isomorphisms along the left edge and
$\gamma(L')$ is given by the composite of the two isomorphisms along
the right edge it follows that $\gamma(L)=\gamma(L')$.   The fact
that $\gamma$ is a group homomorphism when restricted to pure string links
was observed in Section 2.
\qed

We next show that $\gamma(L)$ is trivial on boundary string links.  When
combined with Theorem \ref{iequivalence} this implies that the Gassner
invariant is trivial on string links that are I-equivalent to boundary links.

The most convenient definition of boundary link in the context of string links
is the following.

\begin{definition}  Let $G(x_1,\cdots,x_n)$ denote the free group on the
generators $x_i$.  A string link $L$ with $n$ components  and complement $X=
D^2\times I-L$ is called a {\it boundary string link} provided there exists a
homomorphism 
$$\phi:\pi_1X\to G(x_1,\cdots x_n)$$
so that $$\phi(\mu_i)=x_i=\phi(\mu_i')$$
for each $i=1,\cdots n$ (with the $\mu_i$ as above).
\end{definition}

The standard transversality argument shows that this is equivalent to the
existence of $n$ disjointly embedded surfaces $E_1,\cdots , E_n$ in $D^2\times
I$ so that 
the boundary of each $E_i$ is a circle decomposed into two arcs, one of
which  is the $i$th component of the string link and the other which made up
of three pieces in $\partial(D^2 \times I)$ in the form $A \times \{0\} \cup
\{b\} \times I \cup A \times \{1\}$ where $A \subset D^2$ is an arc from $p_i$
to
$b \in \partial D^2$.

\begin{theorem}
If $L$ is a boundary string link, then $\gamma(L)$ is the identity
transformation.
\end{theorem}
\noindent{\sl Proof.}
Write $G(x_1,\cdots, x_n)=\pi_1(\vee_{i=1}^n S^1)$.  There exists a homotopy 
commutative diagram
\[
 \begin{diagram}
 \node[2]{X_0}\arrow{s,l}{i_0}\arrow{se,t}{\Phi_{|X_0}}\\
\node{\vee_{i=1}^n S^1}\arrow{ne,t}{j_0}\arrow{se,b}{j_1}
\node{X}\arrow{e,t}{\Phi}\node {\vee_{i=1}^n S^1}\\
\node[2]{X_1}\arrow{n,l}{i_1}\arrow{ne,b}{\Phi_{|X_1}} 
 \end{diagram}
\]
where  $\Phi$ induces $\phi$, and so that $j_0(x_i)=\mu_i$ and
$j_1(x_i)=\mu_i'$.  

By definition, the Gassner invariant is the composite 
\begin{equation}\label{oivhjew}
(j_0^*)^{-1}\circ (i_0^*)^{-1}\circ i_1^* \circ
j_1^*
\end{equation}
 where the $*$ superscripts refer to cohomology rel $p$ or $P$
with
$F$ coefficients. But since the diagram commutes up to homotopy, and
$\Phi_{|X_i}$ and $ j_i$ are homotopy inverses for $   i=0 $ and $1$, the
composite of Equation \ref{oivhjew} is equal to the identity, as desired. \qed

The same analysis applies to the Burau representation.  We leave the
formulation to the reader.   

At this point we note an observation of Le Dimet and X.-S. Lin that the braid group
injects into the string link concordance group.  A sketch of the argument  follows. 
Suppose a braid is concordant to the trivial braid.  The fundamental group of the
complement of the braid includes in the group of the complement of the
concordance.  By Stallings' Theorem \cite{stallings} the inclusion induces
isomorphisms on the lower central series of these groups.  Hence, the
automorphism induced by the braid is trivial on the lower central series of the
free group since it coincides with the automorphism induced by the trivial braid. 
Since the free group is residually nilpotent, this automorphism is trivial and hence
the braid is trivial since the automorphism determines the braid.

\section{String link closure and a factorization of the Alexander
polynomial}

In this section we use the results of Section \ref{wekvoi} to derive a
factorization of the Alexander polynomial   of the closure of a string
link  in terms of the   Gassner matrix and a
certain relative Reidemeister torsion invariant. We also give a  general
method to relate the Alexander polynomial of the link
closure and the knot closure of a string link.  (This is closely related to
the factorization found by Levine in 
\cite{levine}.)  The results of this section are consequences of two basic facts. 
The first is that the Alexander matrix of the closure of a string link factors  over
$F$ as a product of two matrices one of which is constructed from the Gassner matrix
of the string link.  The second fact is that the Gassner matrix has  the same
linear algebraic properties as the Alexander matrix of a (closed) link, but is a
well defined concordance invariant.

 In the case of a (pure)
braid, we recover the well known fact that the Alexander polynomial of
the closure of a  braid is equal to $\det(I-\gamma(L)(11))/(1-t_1)$.  
(Here $M(1,1)$ represents the $(1,1)$ minor of the square matrix $M$.) 
Another special (and well-known) case occurs for a
$1$-component string link when our factorization   says that the Alexander
polynomial of the closure (a knot) is equal to a certain  Reidemeister torsion
associated to  the string link. 
 The factorization of the Alexander matrix of the closure also sheds 
light on what concordance information is carried by  the higher
Alexander ideals of the closure.

In the following, given an $n$-component string link $L$ let $\hat{L}$ 
denote the closure of $L$. Thus $\hat{L}$ is the link of $n$  
circles in $\rr^3$ obtained by adding $n$ ``parallel'' strands  in
$\rr^3- D^2\times I$, the $i$th strand  having endpoints $p_i\times
\{0\}$ and $p_i\times \{1\}$, as indicated in Figure 3.

\vskip5ex
\begin{figure}
\centerline{\epsfxsize=1.7in\epsfbox{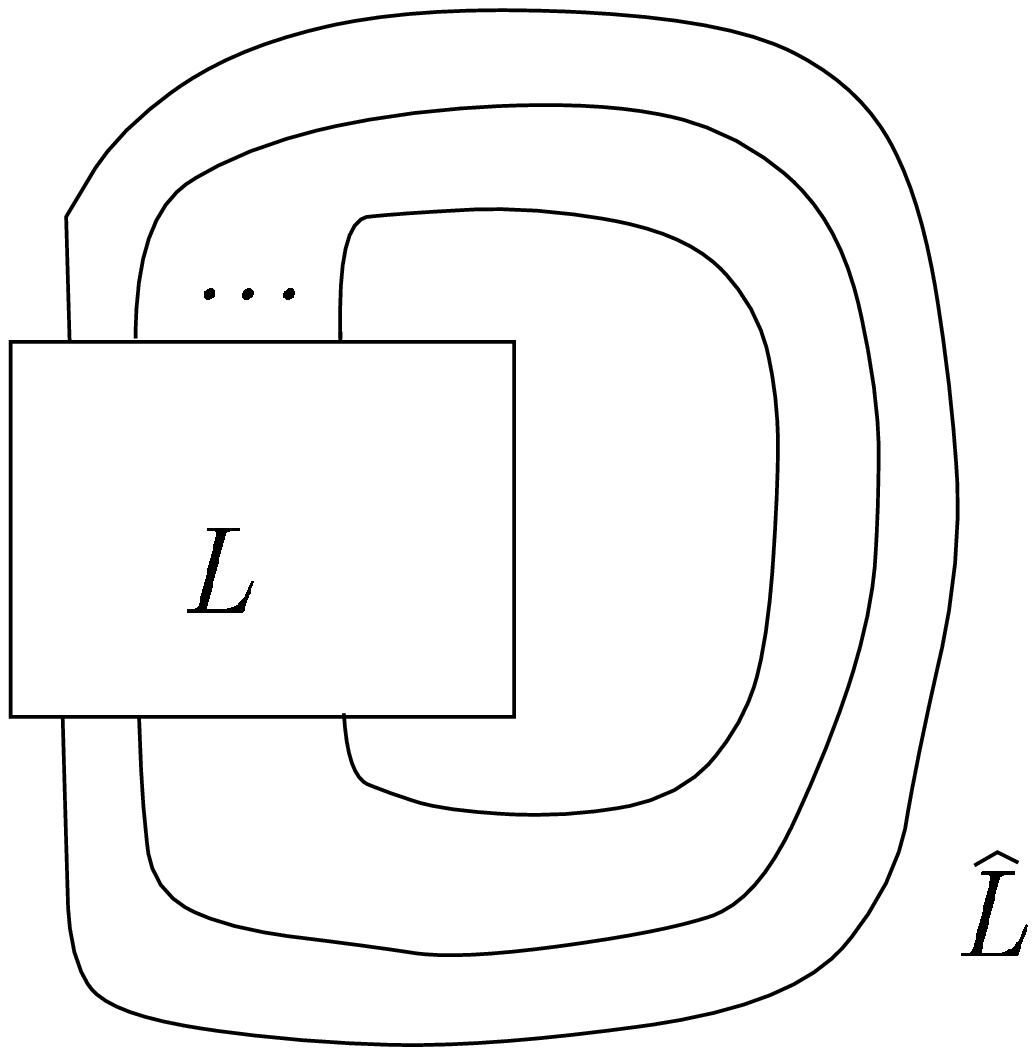}}
\caption{}
\end{figure}
\vskip5ex

We recall briefly the definition of the multi-variable Alexander polynomial of
a link.  Given an $n$-component link ($n\ge 2$) $\hat{L}$ in $S^3$ and a
numbering of its components, let
$\pi=\pi_1(S^3-\hat{L})$ be the fundamental group of the link complement, and
let \break
$\epsilon:\pi\to
\zz^n=<t_1^{\pm1},\cdots,t_n^{\pm 1}>$ be the abelianization map, taking a
meridian to the component colored $i$ to $t_i$.  A presentation of $\pi$
determines, via the Fox calculus, a   matrix with $\zz[\zz^n]$ entries called
the Alexander matrix for the presentation.  The ideal in $\zz[\zz^n]$
generated by the
$(n-1)\times (n-1)$ minors of  this matrix turns out to be of the form 
$\Delta_{\hat{L}}\cdot  J$, where $\Delta_{\hat{L}}\in \zz[\zz^n]$ and $J=\ker
a:\zz[\zz^n]\to \zz$ is the augmentation ideal.   The element   
$\Delta_{\hat{L}}$ is well defined up to units in $\zz[\zz^n]$, i.e. up to
elements of the form $\pm t_1^{a_1}\cdots t_n^{a_n}$ (see \cite{fox2, fox}
and the  paragraphs following the statement of Theorem 6.2 below). It is called the
{\it Alexander polynomial of the link
$\hat{L}$}.

\begin{lemma}
Let $\pmatrix{A&B&C\cr}$ denote the Wirtinger--Fox matrix of a 
(projection of a) string link $L$, as defined in Equation
\ref{alexander}.  Then 
$$V_{\hat{L}}=\pmatrix{A+C& B\cr}$$  is the Alexander matrix of the 
closure
$\hat{L}$ with respect to the Wirtinger presentation of the corresponding
projection of the closure.
\end{lemma}
\noindent{\sl Proof.}  This follows from the fact that the Wirtinger
presentation for the complement of $\hat{L}$ is obtained by setting
$\mu_i=\mu_i'$ in the Wirtinger presentation of the string link $L$.\qed

\begin{theorem}\label{leixntyb}
The Alexander matrix for the closure $\hat{L}$ factors as
 $$V_{\hat{L}}=\pmatrix{A&B\cr}\pmatrix{I_n-\gamma(L)&0\cr
-Z&I_{c-n}\cr}$$
where $I_n$ and  $I_{c-n}$ are identity matrices, and $Z$ is some
$(c-n)\times n$ matrix with $F$ coefficients depending on the
projection. 
\end{theorem}
\noindent{\sl Proof.}  
We have
 \begin{eqnarray*}
V_{\hat{L}}&=&\pmatrix{A&B\cr}+\pmatrix{C&0\cr}\\
&=&\pmatrix{A&B\cr}( I + \pmatrix{A&B\cr}^{-1}\pmatrix{C&0\cr}\ )\\
&=&\pmatrix{A&B\cr}(I-\pmatrix{\gamma(L)&0\cr Z&0\cr})\\
\end{eqnarray*}
The last equality uses Proposition \ref{xys}.\qed

We now examine more closely how to derive the Alexander polynomial of a link
from the Wirtinger presentation associated to a projection of the link. 
What we will see is that the Gassner matrix of a string link shares many
important properties of the Alexander matrix for the Wirtinger presentation of
a (closed) link.

The Wirtinger presentation associates to a projection of (a closed link)
$ {\hat{L}}$ a presentation of
$\pi_1(S^3 - \hat{L})$ such that each generator $x_i$ is a meridian to the
link, and each crossing in the projection defines a relation of the form
$x_ix_jx_i^{-1}x_k^{-1}$. This gives a presentation with $c$ generators and
$c$ relations, where $c$ denotes the number of crossings in the projection.  A
geometric argument shows that any one of the relations is a consequence of the
other $c-1$.  Letting  $V$ denote the Alexander matrix for this 
presentation, a simple computation using the Fox calculus shows that this
implies there exists a vector
$$u=(u_1,\cdots, u_c)\in
\zz[\zz^n]^c$$  each of whose entries $u_i$ is a unit in $\zz[\zz^n]$ so
that 
 $$u V=0.$$  On the other hand, the interpretation of the Fox matrix as a
differential in the group cohomology chain complex  shows that if $a(x_i)=
t_{\alpha(i)}$ then $$Vw=0,$$ where
$$w=\pmatrix{1-t_{\alpha(1)}\cr\vdots\cr1-t_{\alpha(c)}\cr}.$$
We will use this to show that $\Delta_{\hat{L}}$ is equal to 
${{\det (V(i,j))}\over {1-t_{\alpha(j)}} }$, were $M(i,j)$ denotes the matrix
obtained from a matrix $M$ by removing its $i$th row and $j$th column.   

The following lemma can be proven using linear algebra.  We will prove  it (up to
signs)  in Section 10 as an easy application of Reidemeister torsion.

\begin{lemma}\label{oinvxy}
Let $M$ be an $n\times n$ matrix over a domain, and let
$u=(u_1,\cdots,u_n)$ and $w=(w_1,\cdots,w_n)^T$ be a row and column vector
so that 
$$u M=0\mbox{ \ and \ }Mw=0.$$

Then the element  in the fraction field 
$$(-1)^{i+j} {{\det(M(i,j))}\over{u_iw_j}}$$
is independent of $i$  and $j$, i.e. if $u_i,u_p,w_j,w_q$ are all non-zero,
then $$(-1)^{i+j} {{\det(M(i,j))}\over{u_iw_j}}= 
(-1)^{p+q}{{\det(M(p,q))}\over{u_pw_q}}.$$\qed\end{lemma}

From the discussion preceding Lemma \ref{oinvxy} we see that if $V$ denotes
the Alexander matrix of the link $\hat{L}$, then 
$$\Delta_{\hat{L}}(t_1,\cdots,t_n)= {{\det(V(i,j))}\over{1-t_{\alpha(j)}} }.$$

It turns out that the Gassner matrix $\gamma(L)$ of a string link has the same
property as the Alexander matrix of a (closed) link obtained using
the Wirtinger presentation, as the following proposition shows. 

\vfill \eject

\begin{proposition}\label{sadpoi}
 Let $L$ be a string link and $\gamma(L)$ its Gassner matrix.

\begin{enumerate}
\item 
$$(I-\gamma(L))\pmatrix{1-t_1\cr\vdots\cr 1-t_n\cr}=0.$$

\item   
$$(t_1^{-1}, (t_1t_2)^{-1}, \cdots,
(t_1 t_2 \cdots t_n)^{-1})(I-\gamma(L))=0.$$

\end{enumerate}
\end{proposition}
  Proposition \ref{sadpoi} allows us to make the following 
 definition. 

\begin{definition}
Let $L$ be a string link of 2 or more components, with Gassner matrix
$\gamma(L)$.  The {\it Alexander rational function of $L$} is the rational
function 
$$\Delta_L(t_1,\cdots,t_n)=(-1)^{i+j}(t_1\cdots
t_i){{\ \det(I-\gamma(i,j))}\over{ 1-t_j }}.$$ Lemma \ref{oinvxy} implies that
$\Delta_L\in F$   is well defined and independent of
$i$ and
$j$. Of course, it is an $I$-equivalence invariant of the string
link since $\gamma(L)$ is.

\end{definition}
The following lemma is a useful step in proving Proposition \ref{sadpoi}.

\begin{lemma}\label{sdoipogi}
The Gassner matrix for the braid obtained by applying one full positive twist
to the trivial braid is 
$$ (t_1\cdots t_n)(\mbox{Id}+ \pmatrix{1-t_1\cr\vdots\cr 1-t_n\cr} 
(t_1^{-1},\  (t_1t_2)^{-1},\  \cdots,\ 
(t_1 t_2 \cdots t_n)^{-1})).$$

\end{lemma}
\noindent{\sl Sketch of Proof.}  This is an easy computation using the Fox
calculus and the fact that the automorphism of the free group which
corresponds to this braid is given by conjugation by the element
$(t_1\cdots t_n)^{-1}$. We leave the details to the reader.\qed

\noindent{\sl Proof of Proposition \ref{sadpoi}.}
 \begin{enumerate}

\item   Let $X$ denote the complement of $L$. To prove the first statement, it
suffices to show that there is a 1-cocycle on $\pi_1(X)$ which takes the
values
$1-t_i$ on both $\mu_i$ and $\mu_i'$.  But there is an obvious one---the
coboundary of the constant 0-cochain.  Indeed, if $f=1\in F\cong
C^0(\pi_1(X);F)$, then $(\delta f)\in  C^1(\pi_1(X);F)$ is defined by
$(\delta f)(x)= 1-\epsilon(x)$.  Thus $(\delta f)(\mu_i)= 1-t_i=(\delta
f)(\mu_i')$.  

\item  Notice that conjugating a string link by the braid which is obtained by
applying one full positive twist to the trivial braid yields an isotopic
string link. Using Lemma \ref{sdoipogi} this means that 
\begin{equation}\label{sdljsdsad} (t_1\cdots t_n)(\mbox{Id}+ wu)\gamma(L)=
\gamma(L)(t_1\cdots t_n)(\mbox{Id}+ wu)\end{equation} where
$$w= \pmatrix{1-t_1\cr\vdots\cr 1-t_n\cr}$$
and 
$$u=(t_1^{-1}, (t_1t_2)^{-1}, \cdots,
(t_1 t_2 \cdots t_n)^{-1}).$$

Simplifying Equation \ref{sdljsdsad} and using the fact proven in part 1
that $\gamma(L)w=w$ one obtains 
$$wu\ \gamma(L)= wu.$$
Equating the first row of each side of this equation gives
$$(1-t_1) u \gamma(L) = (1-t_1)u.$$
Since $ 1-t_1 $ is non-zero, the desired result follows.
\qed\end{enumerate}

In the case of the 1-variable Alexander polynomial (and/or 1-colored
string links) and the Burau representation, a slight modification of the
definition is necessary.  Let
$\alpha:\zz[\zz^n]\to \zz[\zz]$ be the ring map taking each $t_i$ to
$t$. Let ${\hat{L}}$ be a closed link and $V_{\hat{L}}$ an Alexander matrix
for $ {\hat{L}}$. First, it is not true in general that
det$(\alpha(V_{\hat{L}})(11)/(1-t)$ is a polynomial.  Instead, one defines
the 1-variable Alexander polynomial of a link $\hat{L}$ to be 
$$\bar{\Delta}_{\hat{L}}(t)=\det(\alpha(V_{\hat{L}})(11)).$$   That this is an
invariant follows just as before from Lemma  \ref{oinvxy} along with the fact that
$\alpha(V_{\hat{L}})$ has eigenvectors $\alpha(u)$ on the left and 
$(1,\cdots,1)^T$ on the right (because $(1-t,\cdots,1-t)^T$ is an
eigenvector). Thus  in this situation one defines the 1-variable
Alexander polynomials  by
$$\bar{\Delta}_{\hat{L}}(t)=   \det(\alpha(V_{\hat{L}})(11))$$
and so it is sensible for us to define  the 1-variable Alexander function
for  a string link $L$ to be  
$$\bar{\Delta}_L(t)=  t \det(I-\beta(L)(11)).$$ 
In these definitions it is unnecessary for the components to be colored or
for the string link to be pure.

From the definitions we have (up to multiples of $t$)
$$\bar{\Delta}_{\hat{L}}= (1-t)\alpha(\Delta_{\hat{L}})$$
for closed links 
and 
$$\bar{\Delta}_L=(1-t)\alpha(\Delta_L)$$
for pure (or colored) string links.

A ``coordinate free'' definition of the Alexander function for string links  can be
given in terms of the reduced Gassner invariant and the reduced Burau invariant.

\begin{lemma} \label{noname}
\begin{enumerate}
\item Let $L$ be a pure string link and 
$\tilde{\gamma}(L):H^1(X_0;F)\to H^1(X_0;F)$ the reduced Gassner invariant. Then
$$\Delta_L= \pm {{\
\det(I-\tilde{\gamma}(L))}\over{
  (t_1t_2\cdots t_n)^{-1}-1 }}.$$
\item
Let $L$ be an arbitrary string link and $\tilde{\beta}(L):H^1(X_0;\qq(t))\to
H^1(X_0;\qq(t))$ its reduced Burau invariant. Then

$$
\bar{\Delta}_L={ {\det(I-\tilde{\beta}(L))}\over{ t^{-1}+t^{-2}+\cdots +t^{-n}}}.
$$
 \end{enumerate}
\end{lemma}

We will prove this lemma in Section 10.
\vskip4ex

  We next
use the factorization of Theorem
\ref{leixntyb} to obtain a    factorization of the Alexander
polynomial of the closure of the link in terms of the Alexander
rational function of the string link.

In order to give a topological interpretation of this factorization, we introduce the 
{\it torsion} of a string link.      We use Reidemeister torsion to avoid some
algebraic complications arising from the fact that the ring $\Lambda$ is not a
P.I.D.,  but also  because torsion has many nice properties.  (See Section 10 for an
in depth discussion, including the definition of Reidemeister torsion.  The basic
reference for torsion and its properties is \cite{milnor}.)

Let $L$ be  a string link; as before denote its  exterior by $X$ and the two 
punctured discs in the boundary of $X$ by $X_0$ and $X_1$.  Lemma 
\ref{asdfg} implies that $H^*(X,X_1;F)=0$, so that the chain complex 
$$C^*(X,X_1;F)=\mbox{Hom}_{\Lambda}(C_*(\tilde{X},\tilde{X_1}),F)$$
is acyclic. Here  $\tilde{X}$ refers to the $\zz^n$ cover of $X$.  This chain
complex is based by taking lifts of a (relative) cell structure  on the pair
$(X,X_1)$.  

\begin{definition}
Define $\tau(L)$, the {\it Reidemeister torsion of a string link $L$}, to be the 
Reidemeister torsion of the based, acyclic cochain complex $C^*(X,X_1;F)$.
\end{definition}

 Thus $\tau(L)\in F^*$ and is well defined up to
multiplication by   units in 
$\Lambda$, i.e. up to multiplication by elements of the form  $\pm
t_1^{a_1}\cdots  t_n^{a_n}$.    It is a smooth
invariant of the string link 
$L$. In particular, it is
independent of the choice of cell structure.  Moreover, standard properties of
the torsion show that  $\tau(L)$ is multiplicative with respect to multiplying
string links, i.e.
$\tau(L_1L_2)=\tau(L_1)\tau(L_2)$.  This also follows by simple linear algebra from
the fact proven below that $\tau(L)=\det(A\ B)$, where $(A\ B\ C)$ denotes the
Wirtinger--Fox matrix.

In the 1-variable case one
makes the analogous definition. We denote by
$\bar{\tau}$ the 1-variable torsion of a string link.

Notice that if $L$ is a {\it braid} then $\tau(L)=1$ since $X$ collapses to
$X_1$.

In the following theorem the formulas are understood to hold only modulo 
units in $\zz[\zz^n]$.

\begin{theorem}\label{factorization}
Let $L$ be a pure  string link and let $\hat{L}$ denote its closure.  
Then the (multi-variable) Alexander polynomial of  $\hat{L}$ is the product
of the torsion of the string link $\tau(L)$ and the Alexander rational
function of the string link $\Delta_L$:
$$\Delta_{\hat{L}}=\tau(L)\cdot \Delta_L.$$
Moreover, given a projection of the string link $L$, $\tau(L)=\det(A\ B)$ 
where 
$(A\ B\ C)$ is the Wirtinger--Fox matrix.  In particular,
$\tau(L)\in\zz[\zz^n]$  and 
$\tau(L)$ augments to $\pm 1$.  

Similarly, for a (not necessarily pure)  string link $L$ the single variable
Alexander  polynomial of  $\hat{L}$ factors as
$$\bar{\Delta}_{\hat{L}}(t)=\bar\tau(L) \cdot \bar{\Delta}_L (t). $$ 
\end{theorem}

\noindent{\sl Proof.} Lemma \ref{paorjvc} says that the matrix $(A\ B)$ has
entries  in $\zz[t_1^{\pm 1},\cdots, t_n^{\pm 1}]$ and has determinant
augmenting to 
$\pm 1$.

We first explain why  $\tau(L)=\det(A\ B)$.  
It is easy to see (and well known)  that the Wirtinger presentation comes from a 
cell structure. In the present context, the exterior $X$ of the string link $L$ has a 
cell structure with one $0$-cell, $2n+c$ $1$-cells labeled  
$\mu_1,\cdots,\mu_n,z_1,\cdots,z_n,\mu_1',\cdots,\mu_n'$,    $c$ $2$-cells 
labeled by the crossings in the projection, and no higher dimensional cells.  The 
cell structure for the subcomplex $X_1$ has the same $0$-cell, $1$-cells 
$\mu_1',\cdots,\mu_n'$, and no higher dimensional cells.

From this one concludes that  the short exact sequence of based cochain complexes 
for the pair $(X,X_1)$, $$0\to C^*(X,X_1;F)\to C^*(X;F)\to C^*(X_1;F)\to 0$$ is 

\[ \begin{diagram}\dgARROWLENGTH=2.7ex
\node[2]{0}\node{0}\node{0}\\
\node{0}\arrow{e}\node{F^c}\arrow{n}\arrow{e,t}{j_2}
          \node{ F^c}\arrow{n}\arrow{e}\node{0}\arrow{e}\arrow{n}\node{0}\\
\node{0}\arrow{e}
   \node{F^{c}}\arrow{n,l}{\delta^1_{(X,X_1)}}\arrow{e,t}{j_1}
      \node{F^{n+c}}\arrow{n,l}{\delta^1_X}\arrow{e,t}{i_1}
           \node{F^n}\arrow{e}\arrow{n}\node{0}\\
\node{0}\arrow{e}
   \node{0}\arrow{n}\arrow{e}
              \node{F}\arrow{n,l}{\delta^0_X}\arrow{e,t}{i_0}
           \node{F}\arrow{e}\arrow{n,l}{\delta^0_{X_1}}\node{0}
\end{diagram}\]

By definition, $\tau(L)$ is the determinant of the differential 
$\delta^1_{(X,X_1)}$.  By construction, the matrix $(A\ B\ C)$ represents the 
differential $\delta^1_X$.  Since the relative complex consists of cochains
vanishing on the $\mu_i'$, it follows that
$\delta^1_{(X,X_1)}$ is represented by the matrix $(A\  B)$, as claimed.

For convenience change the projection by isotoping $L$  to add a
little ``kink'' on the first component near its initial point (and  relabel
the
$z_i$) to obtain a new   Wirtinger relation    
$\mu_1 z_1^{-1}$. Then    no  other relation involves $\mu_1$. Taking this to
be the first Wirtinger relation what is achieved is:
\begin{enumerate}
\item The first row of  $(A\
B)$  has $ 1$ in the first entry, $-1 $ in the 
$n+1$st entry, and zeros elsewhere.  
\item The first column of $(A\ B)$ has $1$ in the 
first entry and zeros elsewhere.  
\end{enumerate}  This procedure does not change the determinant 
of $(A\ B)$ up to sign.

Recall that
$$\Delta_{\hat{L}}={1\over{1-t_1}}\det(V_{\hat{L}}(1,1)).$$

Using elementary row operations of the form 
\begin{enumerate}
\item Add a multiple (by an element of $F$) of  the $i$th row to the $j$th row for 
$2\leq i,j\leq c$ 
\item Interchange the $i$th row with the $j$th row for $2\leq i,j, \leq c$
\end{enumerate}
the matrix $(A\ B)$ can be transformed into a matrix of the form
\begin{equation}\label{shgkwuvn}
\pmatrix{1&0\cdots 0\ -1\ 0\cdots 0\cr  \matrix{0\cr  \vdots\cr 0\cr}&D\cr}
\end{equation}
where $D$ is a diagonal matrix.  Note that $\det(D)=\pm \det(A\ B)$.

The elementary row operations used only operate on the second through $c$th 
rows, and hence there exists a $c\times c$ matrix $E$ with determinant $\pm1$ of 
the form
$$E=\pmatrix{1&0\cr 0&X\cr}$$
with $X$ a $(c-1)\times (c-1)$ matrix so that the matrix of Equation 
\ref{shgkwuvn} is equal to $E \cdot (A\ B)$.

Thus
\begin{eqnarray*}
(1-t_1)\Delta_{\hat{L}} &=&\det(V_{\hat{L}}(11))\\ &=
&\det(X\cdot (V_{\hat{L}}(11)))\\
&=&\det( (E\cdot V_{\hat{L}})(11))\\
&=&\det((E (A\ B)\pmatrix{I-\gamma(L)&0\cr-Z&I\cr})(11)) 
\end{eqnarray*}

The last equality follows from Equation \ref{leixntyb}.

A little matrix arithmetic using the fact that the matrix of Equation 
\ref{shgkwuvn} is equal to $E \cdot (A\ B)$ shows that 
$$ (E (A\ B)\pmatrix{I-\gamma(L)&0\cr-Z&I\cr})(11)=
\pmatrix{D_1(I-\gamma(11))&0\cr Z'&D_2\cr}$$
where $D=\pmatrix{D_1&0\cr0&D_2\cr}$,  $Z'$ is some matrix, and 
$\gamma(L)(11)$ is the matrix obtained from $\gamma(L)$ by dropping 
the first row and column.   Thus 
\begin{eqnarray*}
(1-t_1)\Delta_{\hat{L}} &=&\det\pmatrix{D_1(I-\gamma(11))&0\cr 
Z'&D_2\cr}\\
&=&\det(D)\det(I-\gamma(11))\\
&=&\det(A\ B)\det(I-\gamma(11)).
\end{eqnarray*}
Since det$(A\ B)=\tau(L)$, and $\det(I-\gamma(11))/(1-t_1)=\Delta_L$, the
theorem is proved in the multi-variable case.  A similar  
argument works in the 1-variable case. \qed

  Theorem  \ref{leixntyb} can also be used to get  information about the
higher Alexander ideals of the closure of a string link $L$.  Recall that the 
{\it $k$th Elementary Alexander ideal of a link}, $E_k$, is the ideal in
$\zz[\zz^n]$ generated by   $c-k$ minors of an Alexander matrix $V$
for the link, where $c$ is the number of generators of the presentation
defining $V$.  For example, $E_0=0$ since the Alexander matrix  is
singular, and we saw above that $E_1$ is the product of the augmentation
ideal with the Alexander polynomial $\Delta$.

Similar definitions apply to the 1-variable Alexander ideals.

As an immediate consequence of Theorem  \ref{leixntyb} we
get the following.

\begin{corollary}\label{irncymuo}
For $k\leq n$ the first $k$   Alexander ideals $E_k(\hat{L})$  of the closure
$\hat{L}$ of a pure string link $L$   link vanish  if and only if the
Gassner matrix
$\gamma(L)$ has the eigenvalue $1$ with multiplicity $k$.
Similarly, the first $k$ one-variable Alexander polynomials of the
closure
$\hat{L}$ of a string link $L$   link vanish  if and only if the Burau
matrix 
$\beta(L)$ has the eigenvalue $1$ with multiplicity $k$.

In particular, a pure string link is in the kernel of the Gassner   representation
if and only if the first $n$
  Alexander ideals of its closure vanish.  Similar remarks apply to the Burau
representation.
\end{corollary}\qed

We finish this section by deriving a formula relating the (1-variable) Alexander
polynomials of the link and knot closures of a string link, giving an alternative and
more general  approach to the formulas of \cite{levine}.

As pointed out in \cite{levine} any such formula  permits one to compare the
Alexander polynomials of a closed link to the Alexander polynomial of  a knot
obtained by banding the components together; the choice of bands essentially
gives a representation of the link as the closure of a string link.  

Let $L$ be a pure $n$-component string link and  let $B$ be any $n$-braid. 
The case of most interest is when $B$ is 
the braid in Figure 4, since 
the closure $\widehat{LB}$ is a 
reasonable definition for the ``knot closure'' of
$L$.  Clearly $\widehat{LB}$ is obtained 
from  $\hat{L}$ by banding the components together.

\vskip5ex
\begin{figure}
\centerline{\epsfxsize=1.3in\epsfbox{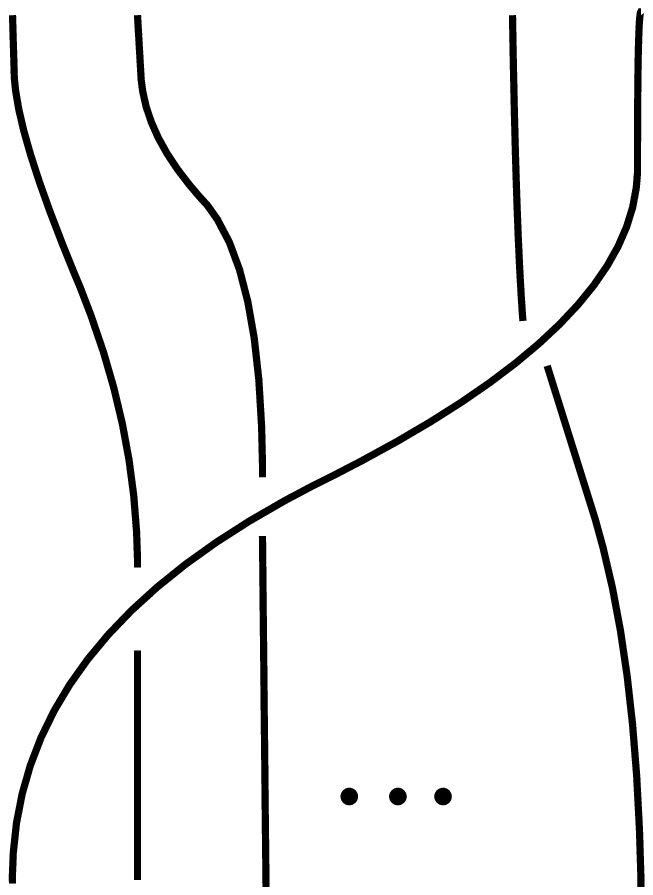}}
\caption{}
\end{figure}
\vskip5ex
Recall that $\tilde{\beta}$ denotes the reduced Burau representation,
$\bar{\Delta}_{\hat{L}}$ denotes the 1-variable Alexander polynomial of the
closed link $\hat{L}$.

\begin{theorem} \label{alexp}
With notation as above, The 1-variable Alexander polynomials of the link closure
$\hat{L}$ and the knot closure $\widehat{LB}$ of a  string link are related by the
formula
$$ \bar{\Delta}_{\hat{L}}= \bar{\Delta}_{\widehat{LB}}\cdot {{\det(I-\tilde{\beta}(L)
)}\over { {\det(I-\tilde{\beta}(L)\tilde{\beta}(B)
)}}}.$$
The correction term $ \det(I-\tilde{\beta}(L)
)/\det(I-\tilde{\beta}(L)\tilde{\beta}(B)
)$ is a string link I-equivalence invariant of $L$.
\end{theorem}
(We use  the reduced Burau representation for convenience.  Alternatively one
can use $\det(I-\beta(11))$  in this formula.)

\vskip3ex
 
\noindent{\sl Proof.}   Theorem \ref{factorization} (in the 1-variable case)  says
that 
$ \bar{\Delta}_{\hat{L}}=\bar{\tau}_L \bar{\Delta}_L $, and similarly for $LB$.  
Since $B$ is a braid, $\bar{\tau}_B=1$ and so
$\bar{\tau}_{LB}=\bar{\tau}_L\bar{\tau}_B=\bar{\tau}_L$.  Moreover, 
${{\det(I-\tilde{\beta}(L))}\over{ ( t^{-1}+\cdots+t^{-n})}}=
\bar{\Delta}_L$ by Lemma \ref{noname}.

Hence
$$
{  {\bar{\Delta}_{\widehat{LB}}    }    \over {   \bar{\Delta}_{\hat{L}}  } }  
= { {\bar{\Delta}_{ {LB}} }\over {\bar{\Delta}_{ L} } }  
={{ \det(I-\tilde{\beta}(LB))}\over { 
\det(I-\tilde{\beta}(L))}}.  
$$
The last assertion follows from Theorem \ref{iequivalence}.  \qed

 \section{Some examples}

We will give a few examples to illustrate how to use the invariants and at
the same time we introduce some geometric constructions on string links.  We will
omit most computations, but remark that they are easily carried out using a
computer algebra package using Equation \ref{sdajasd}.  Given a string link, its
concordance inverse is the string link given by reflecting the given string link
through the plane $\rr^2\times {1\over 2}$ and reversing the orientation of the
components.

A simple argument shows that any string link whose closure is the
trivial link is string link concordant to the trivial link.  Thus if
$L$ is any string link whose closure is trivial, then $\gamma(L)$ is the
identity matrix.

Now consider the pair of links $L_1$ and $L_2$ illustrated in Figure
5.

\vfill
\eject

\vskip5ex
\begin{figure}
\centerline{\epsfxsize=1.6in\epsfbox{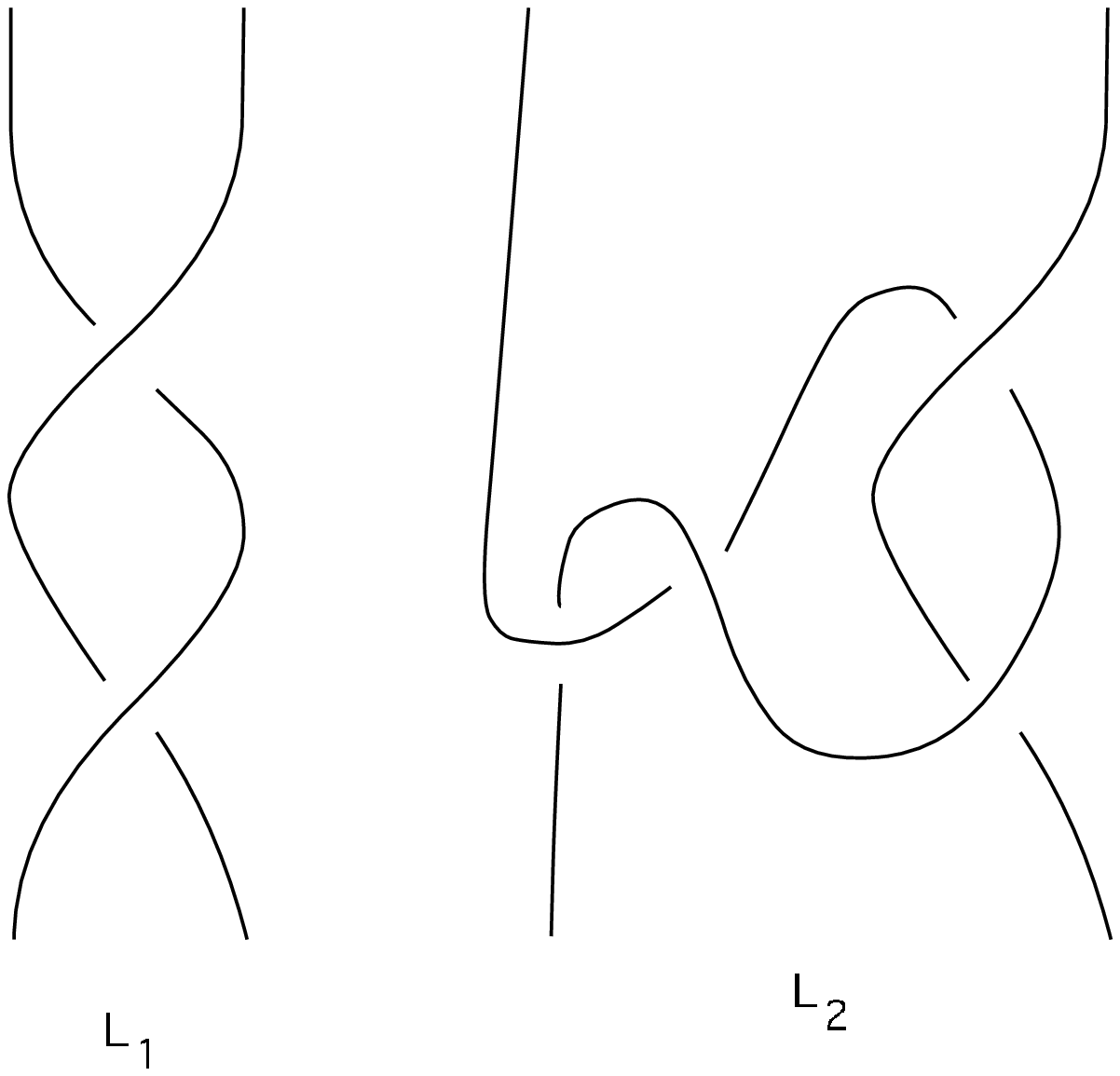}}
\caption{}
\end{figure}
\vskip5ex

Notice that $L_1$ and $L_2$ have the same closure, namely the Hopf link.  
The Gassner matrices for these links are: 
$$\gamma(L_1)=\pmatrix{t_2&1-t_1\cr(1-t_2)t_2&1-t_2+t_2t_1\cr}$$
and 
$$\gamma(L_2)= {1\over {-t_1 - t_2 + t_2 t_1 }}
\pmatrix{  1 - 2 t_2 - t_1 + t_1t_2&t_1-1\cr (t_2-1)t_2&-t_1\cr}$$
so that $L_1$ and $L_2$ are not concordant as string links, even though
their closures are isotopic.

For
a linking number zero example, consider the links in Figure 6.  
\vskip5ex
\begin{figure}
\centerline{\epsfxsize=1.8in\epsfbox{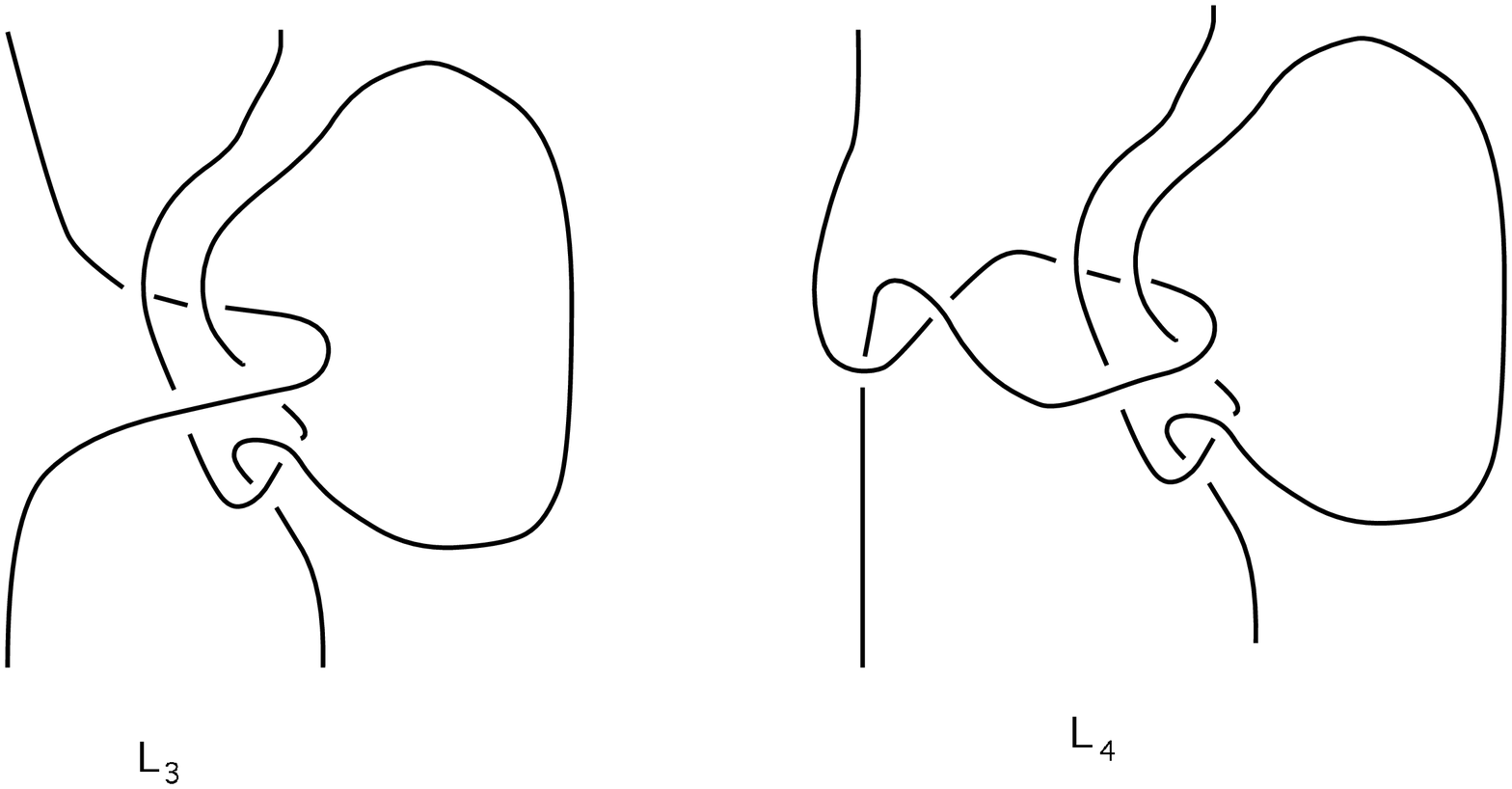}}
\caption{}
\end{figure}
\vskip5ex

The Gassner matrices are complicated but we can distinguish $\gamma(L_3)$ from
$\gamma(L_4)$  using eigenvalues. As with all string links, both
 have $1$ as eigenvalue, the other eigenvalue for $\gamma(L_3)$ (or equivalently
$\tilde{\gamma}(L_3)$) is 
 $${{t_1t_2(2-t_1-t_2+t_1t_2)}\over{1-t_1-t_2+2t_1t_2 }}$$
and the other eigenvalue for $\gamma(L_4)$ is 
$$-{{ -2 t_2t_1-t_1^2+t_1 t_2^2+2 t_1-t_2^2+t_2t_1^2+2t_2-1 }\over
{t_1^2-2t_2t_1^2-t_2-t_1+t_2^2+2t_2t_1-2t_1t_2^2+t_2^2t_1^2}}.$$
Thus  these are not conjugate in the string link concordance group.  (Notice
that for 2-component string links the Gassner representation is abelian.)

These  examples illustrate two interesting constructions 
which can be  used to alter string links.  The first adds a
horizontal twist to one component, as illustrated in going from $L_1$ to $L_2$ in
Figure 5 and in going from
$L_3$ to $L_4$ in Figure 6.  This yields different string links with the same
closure.  More generally  one can add more twists and also add twists to different
components, yielding hosts of examples of string links with the same closure
but different Gassner matrices.  We will show in the next section that in most
cases this construction changes the Gassner matrix.

The second construction is a form of ``Whitehead doubling'' of one component,
as  illustrated in going from
$L_1$ to
$L_3$ and
$L_2
$ to $ L_4$ in Figures 5 and 6.  This is a   subtle construction.  Notice for
example that if one Whitehead doubles each component, the result is a boundary
string link and hence has trivial Gassner matrix.  
 
As another example consider the two 3-component string links
$L_5$ and $L_6$ in Figure 7.  
The first, $L_5$, is a braid.  The second, $L_6$,  is not. 
 
\vskip5ex
\begin{figure}
\centerline{\epsfxsize=1.7in\epsfbox{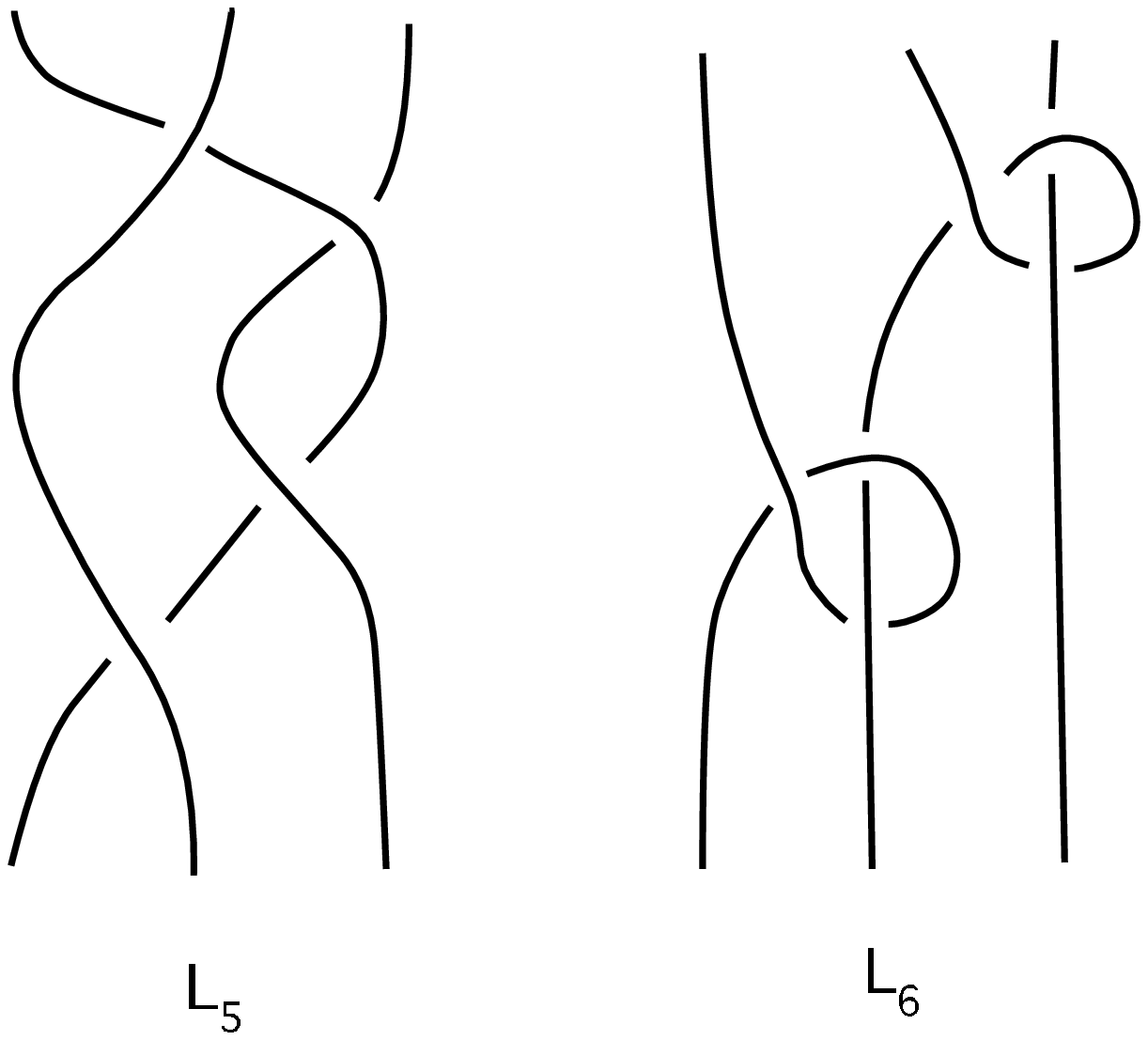}}
\caption{}
\end{figure}
\vskip5ex

A computation shows that the $(1,1)$ entry of $\gamma(L_6 L_5 L_6^{-1})$
is not a Laurent polynomial in   $t_1, t_2, t_3$.  This implies that the
conjugate $L_6 L_5 L_6^{-1}$ of the braid $L_5$ is not string link concordant
to a braid.  Since adding extra ``trivial'' components   replaces the Gassner
matrix by its direct sum with the identity matrix, we   obtain the
following theorem.

\begin{theorem} For $n\ge 3$ The pure braid group on $n$ components is not normal
in the pure string link concordance group.
\qed 
\end{theorem}

(Recall, by the comment at the end of Section 5, the braid group is a subgroup of
the string link concordance group.)

Finally,  we mention an example that shows that in the context of
I-equivalence of string links, the Gassner invariant is not faithful.  In
\cite{orr}  an example of a two component link $\hat{L}$ is given with the property
that it is not concordant to a boundary link even though all its Milnor invariants
are trivial.  A string link, $L$, can be built from this link in such a way that a
direct calculation shows the Gassner matrix is the identity.  On the other hand,
$L$ cannot be I-equivalent to a trivial link.  If it were, then forming the
connected sum of it closure $\hat{L}$ with knots $J_1$ and $J_2$ would yield a link
$\hat{L'}$ that is concordant to the trivial link.  But then forming the connected
sum of
$\hat{L'}$ with the mirror images of $J_1$ and $J_2$ would yield a link $\hat{L''}$
that is concordant to a split link formed as the union of the mirror images of $J_1$
and
$J_2$.  This is a boundary link but is also concordant to $\hat{L}$

\section{Random walks and labelings of string links} In this section we
begin by recalling  a ``probabilistic'' interpretation of the Burau
representation for string links studied in \cite{lin-wang}, and then extend it
to give a similar interpretation of the Gassner representation.  In that
article the authors assign to a string link  diagram with
$n$ strands an
$n\times n$ matrix with coefficients in the field $\qq(t)$ of rational functions
in one indeterminate $t$. The $(i,j)-$th entry of this matrix is a sum
(possibly infinite) over all the paths in the diagram of a string link
starting at
$p_{i}\times \{0\}$ and ending at $p_{j}\times \{1\}$ weighted in the
following way. A path in the diagram is a ``random walk" along the
components (moving in the direction of the strands) of the string link such
that whenever one comes to a crossing, one can either jump up or down, or one
can continue without jumping. Notice there may be arbitrarily long paths; in
the string link of Figure 2 there are paths which loop
around as many times as one likes by jumping up at the right crossing, or going
around the kinks.

The weight $w(p)\in \qq(t)$ of a path $p$ is defined to be the product over
each crossing $c$ along $p$ of $w_{c}$ $$ w(p)=\prod_{c} w_{c},$$ where
$$w_{c}=
\cases{ t^{\pm 1} &  if the path $p$ goes under at the crossing
c \cr  
1-t^{\pm 1} &  if the path $p$ jumps up at the crossing c \cr
 1 &
 if the path $p$ goes over at the crossing c \cr
 0 &  if the path
$p$ jumps down at the crossing c \cr}  $$ 
The sign $\pm 1$ is chosen
to be the sign of the crossing.  Define a matrix by 
$\beta_{i,j} = \sum_p w(p)$ where the sum is over all paths $p$ that begin
at $(p_i ,0)$ and end at $(p_j,1)$.  In \cite{lin-wang}, is shown that this
sum  converges to a rational function for any string link, and that the
resulting matrix is invariant under Reidemeister moves. When restricted to
braids this gives the Burau representation. We will outline simple proofs and
generalize these facts to the Gassner representation below.

The relationship to the approach of the previous sections  is as
follows. Consider the following labeling  scheme
for string link diagrams. Think of a diagram of a string link as a graph in
$R^{2}$ such that every vertex has valence 4 except for the endpoints.  Define
a {\it labeling of the string link with value in a vector space $V$} to be an
assignment of an element of $V$ to each edge in this graph. 

We introduce
a labeling that gives the Gassner representation. For this we must keep track
of which strand one is on, i.e. we work with colored string links. The vector space
$V$ is taken to be the field
$F=\qq(t_{1},\cdots , t_{n})$ of rational functions in $n$ variables. We will
use $t_{o}$ (resp. $t_{u}$) to denote the variable
  with subscript the color of
the over (resp. under) strand in the diagram. 

\begin{theorem}\label{zxlivcuwerhj} Fix a colored string link diagram.

1. Given any collection of $n$ elements $\lambda_{i}'\in F,i=1,\cdots ,n$,
there exists a unique labeling of the diagram with values in $F$ such that 

(a) The labeling at the $i$th top edge (the edge containing the point
$p_{i}\times \{1\}$) is $\lambda_{i}'$.

(b) At a positive  crossing, if  the four edges  are
labeled
$S,E,W,N$ in such a way that the strands of the string link go as the 
diagram (Figure 8),
\vskip5ex
\begin{figure}
\centerline{\epsfxsize=1.5in\epsfbox{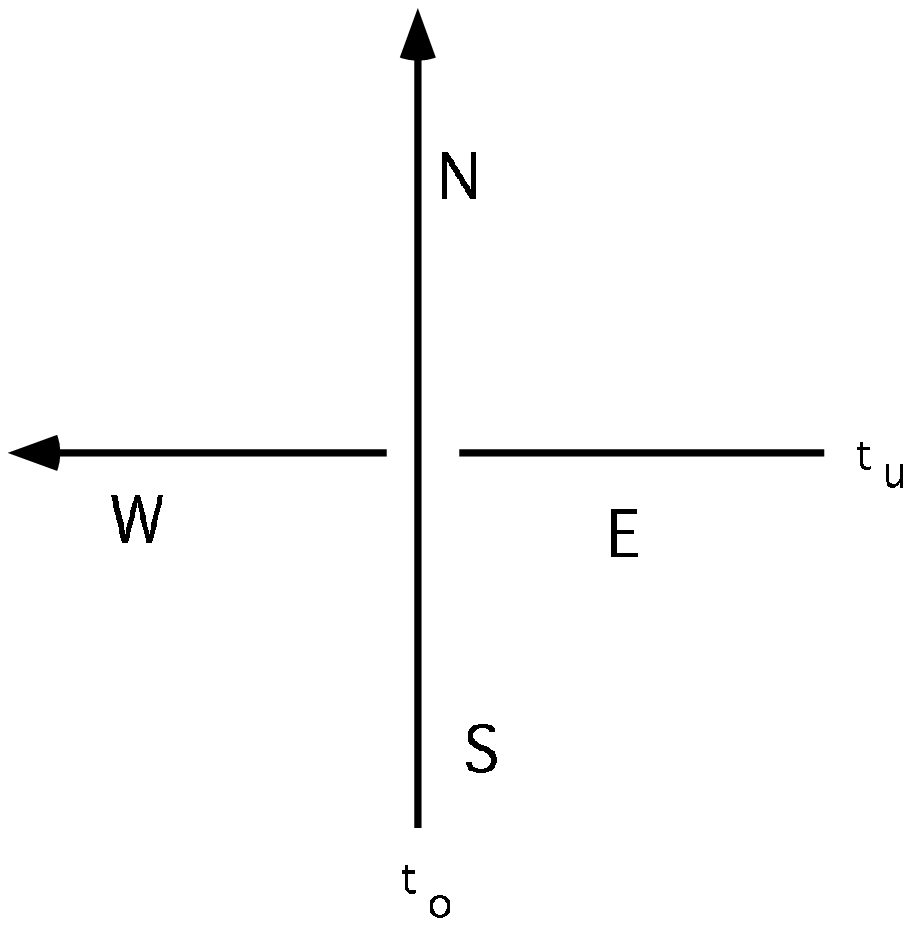}}
\caption{}
\end{figure}
\vskip5ex

then
$S,E,W,N \in F$ satisfy the rules:
 
\begin{eqnarray*} E&=& t_{o} W+(1-t_{u})N \\ S&=&N
\end{eqnarray*}
 
\noindent and at a negative crossing, if the four  edges are labeled 
$S,E,W,N$  in such a way that the strands of the string link go as the
diagram (Figure 9),
\vskip5ex
\begin{figure}
\centerline{\epsfxsize=1.5in\epsfbox{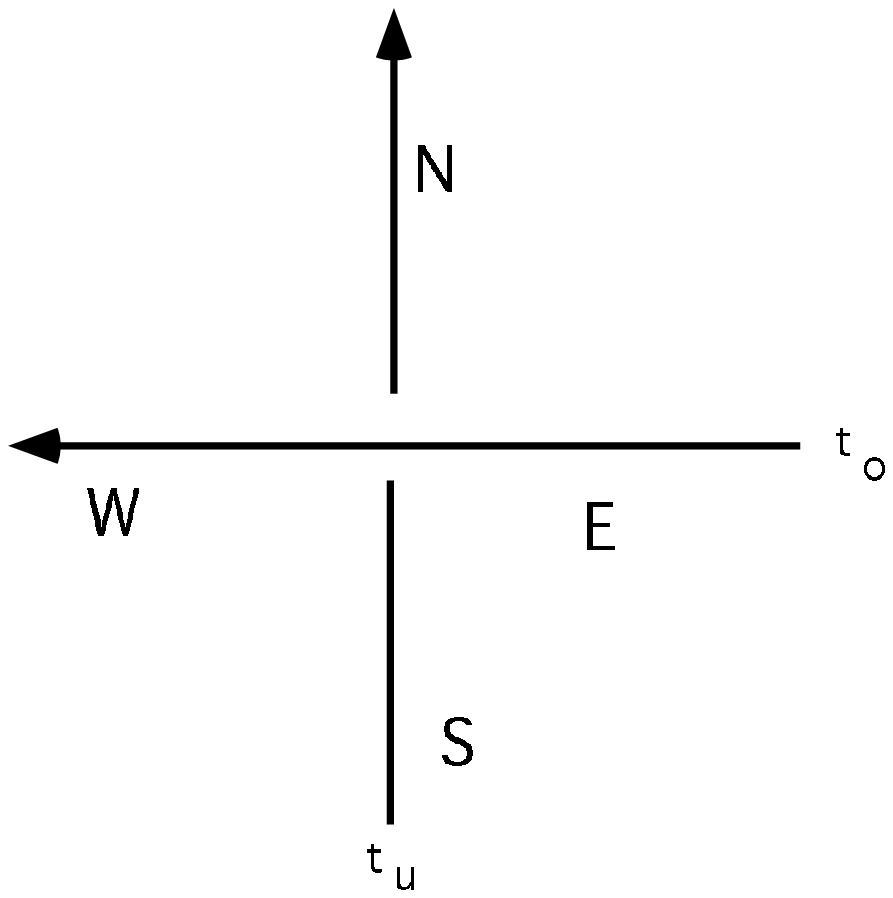}}
\caption{}
\end{figure}
\vskip5ex
 
then
$S,E,W,N$ satisfy the rules:
\begin{eqnarray*} S&=& t_{o}^{-1} N+t_{o}^{-1}t_{u}(1-t^{-1}_{u})W \\ E&=&W.
\end{eqnarray*}

2. Let $G(L)$ be the matrix whose $(i,j)$th entry is equal to the
labeling of the $i$th bottom edge for the unique labeling satisfying the
above equations with $\lambda_{k}'=\delta_{jk}$. Then $G(L)=\gamma(L)$.
\end{theorem}

In terms of random walks, this says that if one assigns the local weight at a
crossing by the rule:
$$w_{c}=
\cases{ t_{o}^{\pm 1} &  if the path $p$ goes under at the
crossing c\cr
  \sigma_{\pm}(1-t_{u}^{\pm 1})& if the path $p$ jumps up at
the crossing c\cr
  1 &  if the path $p$ goes over at the crossing c\cr 
0&  if the path $p$ jumps down at the crossing c \cr} 
$$ 
where the sign is the sign of the crossing and $\sigma_{+}=1,
\sigma_{-}=t_{o}^{-1}t_{u}$, then the sum over all paths starting from the
middle of the edge to the top $i$-th edge gives a labeling where the top
$k$-th edge is labeled by $\delta_{ik}$. 

\vskip2ex

  \noindent{\sl Sketch of Proof of  Theorem \ref{zxlivcuwerhj}. }   The proof is 
virtually the same as that of Lemma \ref{paslapsl}.  Notice that a labeling in the
theorem has the property that the labels on the two edges which pass over a
crossing are equal. Thus there is a one-to-one correspondence
between such labelings of a string link diagram with values in $F$ and
$Z^{1}(\pi, F)$.  Moreover, if we set the $t_i$ equal to $1$ the equations
reduce to $E=W$ and $N=S$ which has a unique solution given a labeling of the
top strands.\qed

The labeling scheme is so flexible that we were led to wonder if there are
other invariants satisfying such a simple linear system. Specifically,  given
any field,   one can seek
  an assignment of local weights at a crossing so that the corresponding
labeling obtained by summing over all paths as above is invariant under
Reidemeister moves.

In general, suppose that some   local weights are assigned as follows: 
$$w_{c}=
\cases{ a_{o,u}^{\pm} & if the path $p$ goes under at the
crossing c\cr
 b_{o,u}^{\pm} &  if the path $p$ jumps up at the crossing
c\cr
 c_{o,u}^{\pm} &  if the path $p$ goes over at the crossing c\cr
d_{o,u}^{\pm} & if the path $p$ jumps down at the crossing c\cr}
  $$ where
$a_{o,u}^{\pm},b_{o,u}^{\pm}, c_{o,u}^{\pm},d_{o,u}^{\pm}$ are arbitrary
elements in some fixed field, and the $\pm$ is determined by the sign of the
crossing c. 

If there exists a string link invariant which assigns to a string link diagram 
a matrix whose $ij$-entry is the sum of the products of weights over all paths
from the
$i$th initial point to the $j$th endpoint of the string link, then the
Reidemeister moves force some complicated algebraic relationships to hold
between the elements
$a_{o,u}^{\pm},b_{o,u}^{\pm}, c_{o,u}^{\pm},d_{o,u}^{\pm}$ of the field. 
 Even if one can find elements satisfying all these algebraic relationships,
there may not exist a solution or it may not be unique if one exists.  
 For these reasons we restricted our
search to choices of weights for which the arguments of Section 4 for the
Gassner representation extend. For simplicity we only considered  choices of weights for which the resulting
link invariants are ``balanced'' in the sense that changing the choice of coloring of
the link gives an equivalent invariant. Thus we required the following two
properties to hold.
\begin{enumerate}

 \item[($C_{1}$)]  (Homogeneity) if we re-index the components by a permutation,
then the resulting matrix is obtained by permuting indices. 

\item[($C_{2}$)] (Nondegeneracy) There is a ring map from the field to $\cc$
so that  the resulting linear system over $\cc$ has a unique solution.
\end{enumerate}

Under these conditions, a laborious calculation  can be used to solve all
relations imposed by invariance under Reidemeister moves. It turns out that there
are essentially only two solutions.  We omit the  entirely
calculational proof.

\begin{theorem}
Let $\nu$ be a  string link invariant taking $n$-component pure string links to 
$n\times n$ matrices over some field, and suppose that $\nu$ is obtained by  taking
a string link to the matrix whose 
$ij$-entry is the sum of the products of weights over all paths
from the
$i$th initial point to the $j$th endpoint of the string link.  If conditions
$(C_1)$ and $(C_2)$ are satisfied, then $\nu$ is either a homomorphic image of the
Gassner representation, or else is determined entirely by the pairwise linking
numbers of the closure of the string link. \qed
\end{theorem}

We finish this section by giving an application of the random walk point of
view for the Gassner invariant.  The following theorem has a simple proof,
although it is not at all obvious how to prove it starting with either the
homological or the Fox calculus definitions of the Gassner invariant.  This
exhibits one advantage of the ``random walk'' approach to the Gassner
representation, namely difficult linear algebra is simplified by the
use of geometric series.  

Given a string link $L$, let $T(L)$ be the string link obtained from $L$ by 
adding a (negative) horizontal twist to the first strand, as explained in the
previous section and illustrated in Figure 10.\vskip5ex
\begin{figure}
\centerline{\epsfxsize=1.5in\epsfbox{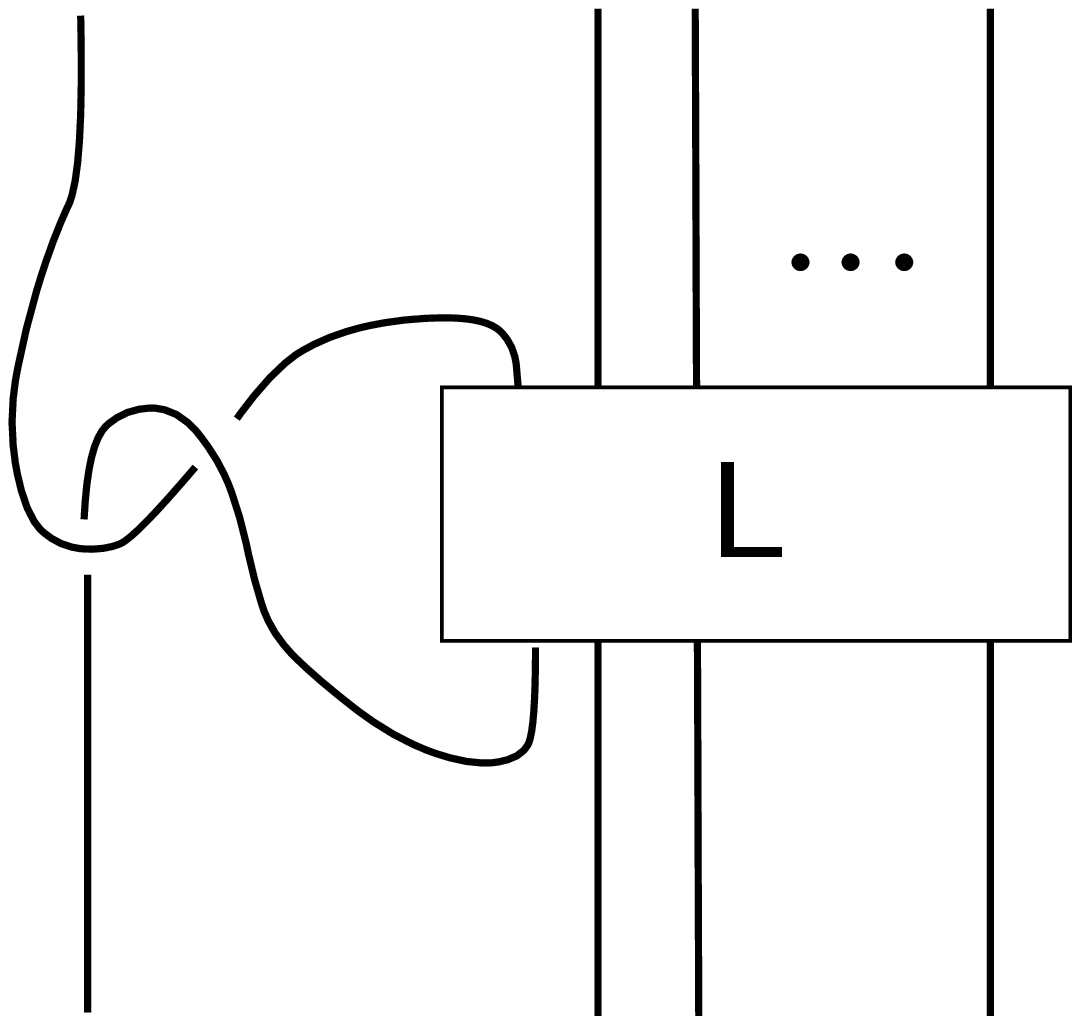}}
\caption{}
\end{figure}
\vskip5ex

The next result shows how to compute the Gassner matrix for $T(L)$ in terms of
the Gassner matrix for $L$.

\begin{theorem}
Let $L$ be an $n$ component pure string link and $T(L)$ obtained from $L$ by
adding a negative horizontal twist to the first strand as illustrated above. 
Suppose that
$$\gamma(L)=\pmatrix{a&{\bf b}\cr \bf{c} &D\cr}$$
where $a\in \Lambda_S$, ${\bf b}$ is an $n-1$ row vector, ${\bf c}$ is an $n-1$
column vector, and $D$ is an $(n-1)\times(n-1)$ matrix. 
Let 
$$\alpha(a)={1\over{1-a(1-t_1^{-1})}}.$$

Then  
$$\gamma(T(L))= \pmatrix{(1-t_1^{-1}) + t_1^{-2}  \alpha(a) a&
t_1^{-1}\alpha(a){\bf b}\cr
 t_1^{-1}\alpha(a){\bf c}& D+(1-t_1^{-1})\alpha(a){\bf cb}\cr}.$$

Moreover, if  $a\ne 1$, then the string links
$L, T(L),T(T(L)),\cdots $ are all distinct in the string link concordance
group but have isotopic closures.  

In general, if $L$ is any pure string link   
with $\gamma(L)\ne I$,  then adding horizontal twists to some component gives an
infinite family of string links with the same closure as $L$ which are pairwise
distinct in the pure string link concordance group.

\end{theorem}

\noindent{\sl Proof.} Let $\gamma_{ij}$ denote the $i,j$ entry of $\gamma(L)$
and let $\gamma'_{ij}$ denote the $i,j$ entry of $\gamma(T(L))$.  Of course
$a=\gamma_{11}$. Call the leftmost crossing   in Figure 10 ``crossing  $1$'',
the next crossing ``crossing 2''; all other crossings of $T(L)$ are crossings
of
$L$ itself.  

 Consider first the entry $\gamma'_{11}$.
We can enumerate all   random walks starting at the initial point of
the first strand and ending at the endpoint of the first strand.  We will
ignore all walks that jump down at any crossing since these have weight zero
and do not contribute to $\gamma'$.  Looking at Figure 10 above, we  see that
there is a path which jumps up at  crossing 1; this path has total
weight $ 1-t_1^{-1} $.  Any other (non-zero weight) random walk passes under  
crossing 1, over   crossing 2, and into the region corresponding to the
string link $L$.   If this random walk is to end up at the endpoint of the
first strand of $T(L)$ it must emerge from the box at the endpoint of the
first strand of $L$.  This walk can then either
\begin{enumerate}
\item  pass under at crossing 2 and
then end, or
\item jump up at crossing 2 and enter the box labeled $L$ again.

\end{enumerate}
 
Thus we can partition the set of all random walks contributing to
$\gamma_{11}'$ according to how many times they pass through the box
labeled $L$.    We already calculated that the (unique) walk never passing
through this box  contributes $ 1-t_1^{-1} $ to $\gamma_{11}'$.  A moment's
thought will convince the reader that the sum of weights over all random walks
passing through the box once is
$t_1^{-1}\gamma_{11}t_1^{-1}=t_1^{-2}\gamma_{11}$.
Similarly the sum of weighs over all random walks
passing through the box twice is
$t_1^{-1}\gamma_{11}(1-t_1^{-1})\gamma_{11}t_1^{-1}=
t_1^{-2}\gamma_{11}^2(1-t_1^{-1})$.
 In general the contribution to $\gamma_{11}'$ of the random walks passing
through the box $n$ times (for $n>0$) is
$t_1^{-2}\gamma_{11}^n(1-t_1^{-1})^{n-1}$.

Summing over $n$ yields
\begin{eqnarray*}
 \gamma_{11}'&=&1-t_1^{-1} + t_1^{-2}\gamma_{11} +\cdots
+t_1^{-2}\gamma_{11}^n(1-t_1^{-1})^{n-1}+\cdots\\
&=&1-t_1^{-1} + t_1^{-2}\gamma_{11}(1+\gamma_{11} (1-t_1^{-1}) +
\gamma_{11}^2(1-t_1^{-1})^{2}+\cdots)\\
&=&1-t_1^{-1} + t_1^{-2}\gamma_{11}\ {1\over {1-\gamma_{11}(1-t_1^{-1})}}\\
&=&1-t_1^{-1} + t_1^{-2}\gamma_{11}\alpha(a)
\end{eqnarray*} 
(The use of geometric series here is what makes the argument using
random walks easy.)

 The other three entries are computed in exactly the same
way.   For example, if $i$ and $j$ are bigger than $1$, then a random walk in
$T(L)$ from the bottom  $i$th endpoint to the top $j$th endpoint can either
pass directly through $L$ to the $j$th endpoint of $L$, contributing
$\gamma_{ij}$, or else it can go to the first endpoint of  $L$ and jump up at
crossing 2 and pass back into the box labeled $L$.  This shows that the
bottom right matrix in
$\gamma(T(L))$ is  the sum of $D$ and some expression involving $\gamma_{i1}$,
$\gamma_{1,j}$, and $\gamma_{11}$.  We leave the verification that
the formula we give is correct to the reader.

Next consider the function 
$$f:\Lambda_S\to \Lambda_S$$
given by the formula 
$$f(a)= 1-t^{-1} + {{t^{-2}a}\over
{1-a(1-t^{-1})}}.$$ Let $z=1-t^{-1}$ and suppose that $a$
has power series expansion 
$$a=a_0 + \sum_{i=N}^\infty a_iz^i$$
where the $a_i\in \qq(t_2,\cdots ,t_n)$ and  $N > 0$.    Since we are assuming $a\ne
1$ we must consider two cases
\begin{enumerate}
\item $a_0\ne 1$ or 
\item $a_0=1$ and $a_N\ne 0$.
\end{enumerate}

A simple calculation shows that 
\begin{equation}\label{function}
f(a) - a= {z\over{1 - za}}  ( 1-a)^2.
\end{equation}
Thus in the first case when $a_0\ne 1$ the power series expansion for $f(a)$ is
given by 
$$f(a)= a_0 + (a_1 + (1-a_0)^2)z + o(z^2)$$
and hence  applying $f$  repeatedly one sees that the $(1,1)$-entry of the Gassner
matrix for $T^k(L)$ is 
$$a_0 + (a_1 + k(1-a_0)^2)z + o(z^2)$$
so that in this case the string links $L, T(L), T(T(L)),\cdots$ are pairwise distinct in
the string link concordance group.

The argument breaks down when $a_0=1$ and so we use a different argument in
this case.  Expanding $f(a)-a$ as a power series using Equation \ref{function} we
find that if
$a_0=1$, the first nontrivial term in the expansion of $f(a) - a$ is $(a_N)^2z^{2N
+1}.$  Hence, the coefficient of
$z^{2N+1}$ in the expansion of
 $f^k(a)=f\circ f\cdots\circ f(a)$ is $a_{2N+1} + k(a_N)^2$. 
In particular these are distinct for distinct values of $k$ and the result is
proved in this case.

The last assertion of the theorem follows from Theorem \ref{traceresult} since one
can  add a twist to any component, not just the first.

\qed

\section{Dominance by finite type invariants} 

If we are allowed to replace some crossings of a string link with transverse
double points, we get a singular string link. Applying the (now standard) method
of extending link invariants to singular link invariants, the Gassner
representation can be extended canonically to singular string links. In this
section we will show that the coefficients in a Taylor expansion of the Gassner
representation are finite type invariants. 

Fix a string link with
$n$ components and let
$z_{i}=1-t_{i}$. For a multiple index $I=(i_{1}\cdots i_{j})$, with $i_{k}\in
\{1,2,\cdots, n\}$, we denote
$|I|=i_{1}+\cdots +i_{j}$, and $z^{I}=z_{i_{1}}\cdots z_{i_{j}}$. Taking the
Taylor expansion of the   matrix  $\gamma(L)$ around $(t_1,t_2,\cdots,
t_n)=(1,1,\cdots, 1)$, we obtain 
$$\gamma(L)=\sum_{|I|=0}^{\infty} b_{I}(L)z^{I}.$$

The coefficients $b_{I}(L)$ are obviously 
invariants of string links. An invariant of string links is said to be  
finite type of order $\leq k$ if it vanishes on any singular string link with
more than $k$ double points. See \cite{barnatan}. 

\begin{theorem}\label{poiewr} The invariant $b_{I}(L)$ is of finite type of order
$\leq |I|$. \end{theorem}

\noindent{\sl Proof.} Suppose we are given a string link diagram for a string link
$L$ and suppose that 
$k$ crossings in this diagram are chosen.   There are $2^k$ string links
obtained by changing the signs of these crossings, labeled by a sequence of
signs $\sigma\in(\pm)^k$. Call these various (embedded) string links
$L_\sigma$.

 For convenience,  add a small kink near the
endpoints of each strand. The proof will follow from a manipulation of an
equation like Equation
\ref{sdajasd}, but, since we wish to relate string links whose diagram differ at
some crossings, it is most convenient to use  presentations of the fundamental
groups of the complements of the $L_\sigma$  with the same set of generators.  
  Thus
 instead of the Wirtinger presentation we will use a larger set of generators for
the fundamental group.  Choose one generator for each edge in the 4-valent graph
given by a (generic) projection of the string link.   As before, label the
bottom strands
$\mu_1,\cdots, \mu_n$ and the top strands $\mu_1',\cdots,\mu_n'$, and call the
intermediate generators say $w_i$.  An example of the labeling scheme is indicated in
Figure 11.  

\vskip5ex
\begin{figure}
\centerline{\epsfxsize=1.5in\epsfbox{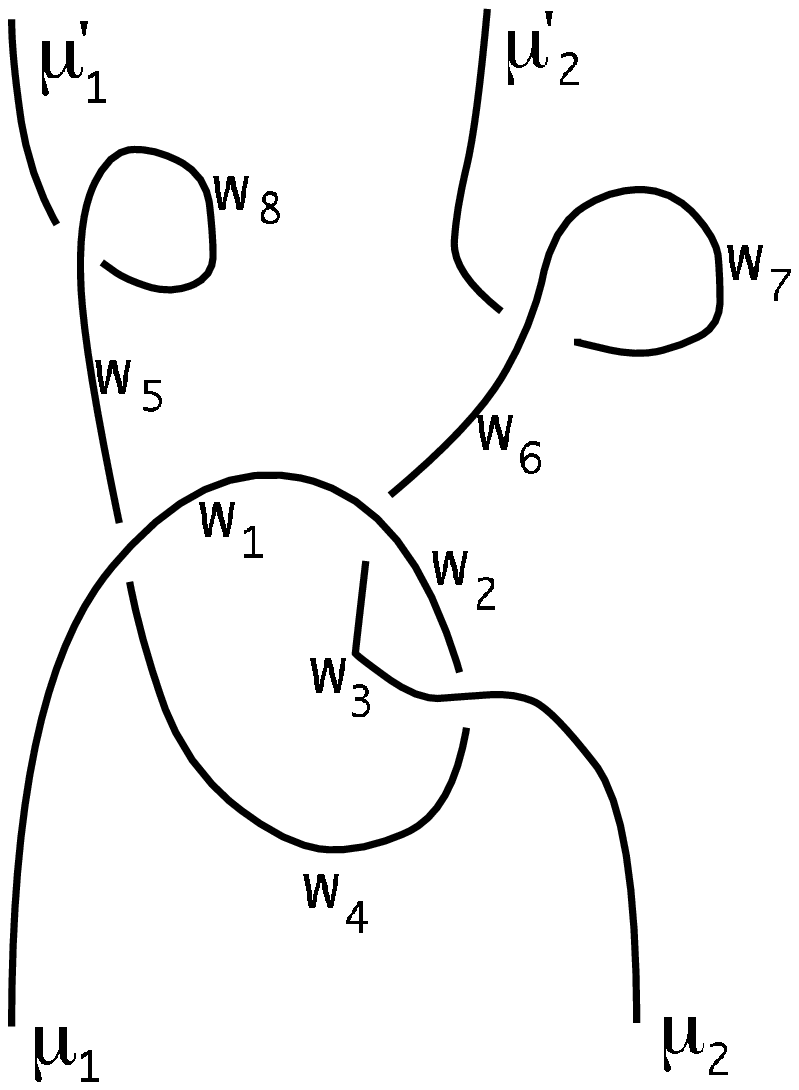}}
\caption{}
\end{figure}
\vskip5ex

For each $\sigma\in (\pm)^k$, the fundamental group of  the complement of  
$L_\sigma$  has a presentation similar to the Wirtinger presentation, but with
extra relations obtained by setting two generators equal if they correspond to
the same ``overstrand''.  Each crossing determines two relations.  In the example
of Figure 11, the first crossing along the first strand determines the
relations 
$\mu_1 w_1^{-1}$ and $w_1 w_5 \mu_1^{-1} w_4^{-1}$.

Applying the Fox calculus to the resulting  presentation for $\pi_1(D^2\times I-
L_\sigma)$ with respect to the ordered basis
$\{ \mu_1,\cdots, \mu_n,w_1,\cdots,w_r,\mu_1',\cdots,\mu_n'\}$ one obtains as
before a matrix
$$\pmatrix{A'_\sigma&B_\sigma'&C_\sigma'\cr}$$
and the identical argument as the one given in Section  4 shows that for each
$\sigma$ there exists a matrix $W_\sigma$ so that
\begin{equation}\label{zxopweajkr}
\pmatrix{\gamma(L_\sigma)\cr
W_\sigma\cr}=-\pmatrix{A_\sigma'&B_\sigma'\cr}^{-1}C_\sigma.
\end{equation}
For notational ease we   rewrite 
$$M_\sigma=\pmatrix{A_\sigma'&B_\sigma'\cr}$$
and 
$$X_\sigma=\pmatrix{\gamma(L_\sigma)\cr
W_\sigma\cr}$$
so that Equation \ref{zxopweajkr} can be rewritten as
\begin{equation}\label{ocghtmsp}
M_\sigma X_\sigma=-C_\sigma'
\end{equation}

\begin{lemma}\label{wemrn}
In Equation \ref{ocghtmsp},
\begin{enumerate}
\item The matrix $C_\sigma'$ is independent of $\sigma\in (\pm)^k$.
\item  The matrix obtained from $M_\sigma$ by mapping each $t_i$ to $1$ has
all entries either $1, -1$ or $0$ and has determinant equal to $\pm 1$. 

\item If $\sigma,\sigma'\in (\pm)^k$ are two sequences of signs that differ
only in the $j$th entry  so that $\sigma_j=1$ and $\sigma'_j=-1$, , then 
$M_\sigma-M_{\sigma'}$ depends only on
$j$ and moreover the Taylor expansion  of $M_\sigma-M_{\sigma'}$ at $t_i=1$ has
zero constant term.
\end{enumerate}
\end{lemma}
\noindent{\sl Proof.}  1.  This follows  from the fact that we added kinks near
the endpoints of each strand, so that the relations involving the $\mu_i'$ are
independent of $\sigma$.

\noindent 2. This follows   from the fact that the pair of relations
determined by each crossing are of the type 
$ab^{-1}$ and $acb^{-1}d^{-1}$ and the same argument   as Lemma \ref{paorjvc}.

\noindent 3.  The matrix $M_\sigma-M_{\sigma'}$ has zeros for all entries except
for a $2\times 4$ submatrix corresponding to the two relations and four
generators involved in the $j$th crossing.  A simple computation (easily derived
from the formulas of Theorem \ref{zxlivcuwerhj}) finishes the argument.
\qed

We continue with the proof of Theorem \ref{poiewr}.  
Let $|\sigma|$ denote the number of
$+$ signs in
$\sigma$. 
Use $o(z^{l})$ to denote any power series with only terms of  order $z^{I}$
for $|I|\geq l$. To prove that
$b_{I}(L)$ is of finite type of order $\leq k$, it suffices to show $$
\sum_{\sigma \in (\pm)^{k}} (-1)^{|\sigma|} X_{\sigma}= o(z^{k}).$$ 

From the first assertion of Lemma \ref{wemrn} we see that 
\begin{equation}\label{wemrnvcoi}
 \sum_{\sigma
\in(\pm)^{k}} (-1)^{|\sigma|} M_{\sigma} X_{\sigma}=0.\end{equation}

We will now prove the theorem by induction.

For $k=1$,
\begin{eqnarray*}
0= M_{+}X_{+}-M_{-}X_{-}&=&(M_{+}-M_{-})X_{+}+M_{-}(X_{+}-X_{-})\\
&=&M_{-}(X_{+}-X_{-})+o(z)
\end{eqnarray*}

 Since $M_{-}=P + o(z)$ where $P$ is an invertible integer matrix, its inverse
is of the form $Q+o(z)$ where $Q$ is the inverse  matrix.  Thus
$$X_{+}-X_{-}=o(z)$$
as desired.   

For the general case, assume the theorem is true for fewer than $k$ crossings.
Expanding Equation \ref{wemrnvcoi} according to the last sign, we have
 $$
\sum_{\sigma' \in (\pm)^{k-1}} (-1)^{|\sigma'|} (M_{\sigma'+}X_{\sigma'+}
-M_{\sigma'-}X_{\sigma'-})=0.$$
This can be rewritten as
$$ 0=\sum_{\sigma' \in (\pm)^{k-1}}
(-1)^{|\sigma'|} M_{\sigma'-}(X_{\sigma'+} -X_{\sigma'-})
+\sum_{\sigma'\in (\pm)^{k-1}} (M_{\sigma'+}-M_{\sigma'-}) (-1)^{|\sigma'|}
X_{\sigma'+}.$$
Now notice that in the second summand, the term $(M_{\sigma'+}-M_{\sigma'-})$ is
independent of $\sigma'$ by the third assertion of Lemma \ref{wemrn}.  Thus the
second summand is $o(z^k):$
\begin{eqnarray*}
\sum_{\sigma'\in (\pm)^{k-1}} (M_{\sigma'+}-M_{\sigma'-}) (-1)^{|\sigma'|}
X_{\sigma'+}&=&(M_{\sigma'+}-M_{\sigma'-})
\sum_{\sigma'\in (\pm)^{k-1}} (-1)^{|\sigma'|}
X_{\sigma'+}\\
&=&(M_{\sigma'+}-M_{\sigma'-})\cdot o(z^{k-1})\ \text{by induction}\\
&=& o(z^{k})\ \text{by the third assertion of Lemma \ref{wemrn}.}
\end{eqnarray*}

Hence
\begin{equation}\label{ewrmn}
o(z^k)=\sum_{\sigma' \in (\pm)^{k-1}}
(-1)^{|\sigma'|} M_{\sigma'-}(X_{\sigma'+} -X_{\sigma'-}) 
\end{equation}

We   continue, inducting on the length of $\sigma$. The goal is to successively
eliminate all occurrences of $M_\sigma$ except for the one sequence
$\sigma=(-,-,-,\cdots,-)$. 

The next step is, using Equation \ref{ewrmn},
\begin{eqnarray*}o(z^k) &=&\sum_{\sigma'' \in (\pm)^{k-2}} (-1)^{|\sigma''|}
M_{\sigma''+-}(X_{\sigma''++} -X_{\sigma''+-}) +M_{\sigma''--}(X_{\sigma''-+}
-X_{\sigma''--})\\
&=&  \sum_{\sigma'' \in (\pm)^{k-2}} (-1)^{|\sigma''|}
M_{\sigma''--}(X_{\sigma''++} -X_{\sigma''+-} +X_{\sigma''-+}
-X_{\sigma''--}) \\
& &\quad +\sum_{\sigma'' \in (\pm)^{k-2}} (-1)^{|\sigma''|}
(M_{\sigma''+-}-M_{\sigma''--})(X_{\sigma''++} -X_{\sigma''+-}) 
\end{eqnarray*}
Using the same argument as before, we see that the second summand in the last
line is
$o(z^k)$, since $(M_{\sigma''+-}-M_{\sigma''--})=o(z)$ and is independent of
$\sigma''$, and because
$$\sum_{\sigma'' \in (\pm)^{k-2}} (-1)^{|\sigma''|}
 (X_{\sigma''++} -X_{\sigma''+-})=o(z^{k-1})$$ by induction (the penultimate
sign is fixed). 

 Thus we obtain
$$
o(z^k)=  \sum_{\sigma'' \in (\pm)^{k-2}} (-1)^{|\sigma''|}
M_{\sigma''--}(X_{\sigma''++} -X_{\sigma''+-} +X_{\sigma''-+}
-X_{\sigma''--})
$$

An induction argument leads to  
\begin{equation} \label{klsadjfrwoieru} o(z^k)=\sum_{\sigma \in
(\pm)^{k}} (-1)^{|\sigma|} M_{-\cdots -}X_{\sigma}. \end{equation}

The inverse of $M_{-\cdots -}$ has a Taylor expansion of the form $Q+o(z)$ for
an integer matrix $Q$, and so inverting $M_{-\cdots -}$ in Equation
\ref{klsadjfrwoieru} finishes the proof.\qed

\section{Torsion and the Alexander function}

In this section we re-examine the results of  
Section 6 using the more sophisticated methods of
Reidemeister torsion.    In particular, the factorization
of the Alexander polynomial of a link will be seen as
an easy consequence of  a Mayer-Vietoris formula for
Reidemeister torsion. Our purpose in including this
material is multifold: for expository reasons, to sharpen
the results of Theorem \ref{factorization}, to obtain a topological
proof, to give a topological interpretation of the determinant of $(A,B)$, and to
fit the Gassner representation for string links in the broader context of knot and
link theory.  We also us torsion to give simple and appealing proofs of Lemmas
\ref{oinvxy} and \ref{noname}.  It is also our belief that the methods of torsion 
are underused in knot theory and often provide a powerful framework to analyze knots
and links.  (Exceptions include Fox and Milnor's use of the interpretation of the
Alexander polynomial in terms of torsion to obstruct knot slicing \cite{fox-milnor}
and
\cite{kirk-livingston} in which a similar interpretation of a twisted Alexander
polynomial related to Casson-Gordon invariants is used to obstruct slicing
algebraically slice knots.  Another important reference is 
\cite{turaev}, in which Reidemeister torsion is applied to the study of links to
vastly simplify much of the previous work on link polynomials as well as to extend
the theory.) 

 For convenience we work with chain complexes and
homology instead of cochain complexes and cohomology.  As remarked in Theorem 3.1
the homology and cohomology Gassner representations are equivalent.

In its most common form, Reidemeister torsion is an
invariant of a based, acyclic complex over a field. Milnor
in \cite{milnor} shows how to extend the definition to  
non-acyclic complexes,  provided one chooses in addition a basis for
the homology of the complex.  From another point of view  the
Reidemeister torsion provides an acceptable substitute for the {\it
order} of a torsion module.  It is
in this full generality that the various facets of Alexander type
invariants of knots and links can be related.

We recall the  definition of $\tau$.  First, if
$V$ is   a vector space over $F$  and ${\bf b}$, ${\bf c}$ are bases
of $V$, let $[{\bf c}|{\bf b}]$ denote the determinant of the
transition matrix 
$M$ where $c_i=\sum_j M_{i,j} b_j$.   Now suppose that $C_*$ is a free
based chain complex over a field. Denote  the basis for $C_i$  by
${\bf e}_i$.  Suppose also that   
bases
$[{\bf h}_i]$ of  $H_i(C_*)$ is chosen.   

We arbitrarily choose bases ${\bf b}_i$ for the boundaries
$B_i=\mbox{Image }\ \partial_{i+1}\subset C_i$.     We also choose
arbitrary lifts
$\tilde{{\bf b}}_i\in C_{i+1}$ of the ${\bf b}_i$.  Finally we
choose representative cycles ${\bf h}_i$ of the classes $[{\bf
h}_i]$. 

With respect these choices the torsion is defined to be the product
$$\prod_i [ {\bf b}_i {\bf h}_i  \tilde{{\bf b}}_{i-1}| {\bf e}_i]^{(-1)^i}.$$ It is
well defined in $F^*/(\pm1)$  (for a
ring
$A$ let
$A^*$ denote the group of units in $A$); in particular it is independent of the choice of the
${\bf b}_i$, their lifts
$\tilde{{\bf b}}_i$, and the choice of cycles, ${\bf h}_i$, representing the
homology bases
$[{\bf h}_i]$.  For a cochain complex, defines the torsion by turning a
cochain complex into a (negatively graded) chain complex via the trick
$C_{i}=C^{-i}$.

Next suppose that $R$ is a domain and let $F$
denote its quotient field. Let $C_*$ be a finitely generated {\it
free } chain complex over
$R$.  Notice that
$H_i(C_*\otimes_{R}F) =H_i(C_*)\otimes_R F $ since $F$ is flat
over $R$.  To deal with the situation when $C_*$ is not
acyclic, suppose that for each $i$ a finitely generated {\it free}
$R$ submodule $K_i\subset H_i(C_*)$ is given such that the inclusion
$K_i\subset H_i(C_*)$ induces an isomorphism $K_i\otimes F\to
H_i(C_*\otimes F)$.   There may be many
choices for $K_i$ and different choices will lead to a different
torsion in general.

Then   choosing bases ${\bf e}_i$ of 
$C_i$ and
${\bf h}_i$ of $K_i$  determines bases (over $F$) of  $C_*\otimes F$ and
$H_i(C_*\otimes F)$. As explained above, these bases determine a
Reidemeister torsion $\tau(C_*\otimes F;{\bf h}_i)\in F^*$.

The  easy but crucial observation is that   the image of   $\tau$
in  
$F^*/R^*$ is independent of the   bases ${\bf e}_i$.   It does
depend on the choice of
$K_i$, but not of the choice of bases $[{\bf h}_i]$ of $K_i$ since changing bases of
$K_i$ (over $R$) gives a determinant in $R^*$.  For convenience denote the image
of
$\tau(C_*\otimes F;{\bf h}_i)$  by
$\tau(C_*;K)\in F^*/R^*$.  If $C_*\otimes F$ is acyclic we
abbreviate this to
$\tau(C_*)$.    Moreover, if $C_*=C_*(S;R)$ for some space $S$ we will use
$\tau(S;K)=\tau(C_*;K)$. 

As an example, suppose that $M$ is a torsion module over $R$, and
suppose that $M$ admits a finitely generated free $R$ resolution
$0\to C_n\to\cdots\to C_0\to M\to 0$.  Then $C_n\otimes F$ is
acyclic and we call $\tau(C_*)$ the {\it order of $M$}, denoted by
$|M|$.  A standard argument shows that
$|M|$ is well defined, i.e. independent of the choice of
resolution.  Notice that if $M$ is a PID then $M$ admits a free
resolution $0\to C_1\mapright{\partial_1} C_0\to M\to 0$.  By
definition $|M|=\det(\partial_1)$ (mod $R^*$), coinciding with the
usual definition of the order of a torsion module.  From this perspective the
Reidemeister torsion provides an extension  to $\zz[\zz^n]$ of the standard knot
theory  methods applied to the P.I.D. $\qq[\zz]$.

Suppose now that 
$$0\to C_*\to D_*\to E_*\to 0$$
is a short exact sequence of free, finitely generated 
$R$-complexes.  Assume free submodules $K_i^C\subset H_i(C_*)$,
$K_i^D\subset H_i(D_*)$, and $K_i^E\subset H_i(E_*)$ are specified
inducing isomorphisms when tensoring with $F$, so that the
torsions $\tau(C_*; K^C)$, $\tau(D_*; K^D)$, and $\tau(E_*; K^E)$ are
defined.   Choosing bases for the $K$ determines bases for the $F$
homology of the complexes, and hence the homology long exact
sequence can be viewed as a  based, acyclic $F$ complex ${\cal H}_*$.
\[
 \begin{diagram}\dgARROWLENGTH=2.3ex
 \node{\cdots}\arrow{e}\node{H_i(C_*\otimes F)}\arrow{e}
\node{H_i(D_*\otimes F)}\arrow{e}\node{H_i(E_*\otimes
F)}\arrow{e}\node{\cdots}\\
\node[2]{K_i^C\otimes F}\arrow{n,r}{\cong}\node{K_i^D\otimes
F}\arrow{n,r}{\cong}\node{K_i^E\otimes F}\arrow{n,r}{\cong}
 \end{diagram}
\]
(The grading is defined by    the convention ${\cal H}_{3i}=H_i(E_*), {\cal
H}_{3i+1}=H_i(D_*)$, and ${\cal H}_{3i+2}=H_i(C_*)$.) We denote its torsion by
$\tau({\cal H})$. It is well defined, independent of the choice of  bases for the $K$.
The same argument as in \cite[Theorem 3.2]{milnor} shows that 
\begin{equation}\label{torsionformula} 
 {\tau(C_*;K^C)\tau(E_*;K^E)} \tau({\cal H})=
{\tau(D_*;K^D)}.
\end{equation}

We will also use a slight extension of the preceding formula.  If $0\to C_*\to
D_*\to E_*\to 0$ is a short exact sequence of  {\it based} chain  complexes {\it
over a field}, and homology bases are chosen, then the Formula
\ref{torsionformula} is still true (and holds in $F^*/\pm1$) provided the bases are
chosen {\it compatibly}.   What this means is that  if ${\bf c, \ d,\ e}$ are the given 
bases for
$C_*,\ D_*,\ E_*$ respectively and ${\bf \tilde{e}}\subset D_*$ is a lift of ${\bf e}$
to $D_*$ then one requires that the determinant of the change of basis matrix from
${\bf d}$ to
$\{ {\bf c,
\tilde{e}}\}$ should be $\pm 1$.

\vskip3ex

We will now apply these ideas to the situation of string links. 
First we prove Lemmas \ref{oinvxy} and \ref{noname}.  These provide a gentle
introduction to computing with torsion.

\vskip3ex

\noindent{\sl Proof of Lemma \ref{oinvxy}.}  Assume that $M$ has rank
$n-1$ or else the lemma is trivial. Consider the following acyclic based complex:
\begin{equation}\label{prejmv}
0\to
F\mapright{W}F^{n}\mapright{M}F^n\mapright{q}F\to 0
\end{equation} where $F$ denotes the
fraction field, $W(f)= f w$, and $q(f_1,\cdots f_n)=
\sum_i u_i f_i$.  We base the complex by taking the standard basis ${\bf e}=\{e_i\}$
in the vector spaces  $F^n$ and the unit $1$ as a basis for $F$.   

To compute the torsion we need to choose bases for the boundaries and lifts
thereof.    Write the complex (\ref{prejmv}) as  
\begin{equation}\label{prejmvi}
0\to D_3\to D_2\to D_1\to D_0\to 0.
\end{equation}
 Then let ${\bf \tilde{b}_2}=1$ so ${\bf b}_2=w$. Letting ${\bf \tilde{b}_1}=
\{e_1,\cdots , \hat{e}_j,\cdots , e_n\}$ one has that ${\bf b}_1=
\{Me_1,\cdots , \widehat{Me}_j,\cdots , Me_n\}$. Let $ {\bf \tilde{b}_0}=  e_i$
so that
${\bf b}_0= u_i$.  Then 
\begin{eqnarray*}\tau(D_*)&=&{{[u_i|1]\ [{\bf b}_2 {\bf \tilde{b}}_1| {\bf e}]}
\over
{[{\bf b_1}{\bf \tilde{b}}_0}|{\bf e}] [1|1]    }\\
&=& {{u_i\cdot  (-1)^{j+1}w_j}\over {(-1)^{i+n}\det(M(ij))\cdot 1}}.
\end{eqnarray*}
Since the left side is well defined and  independent of $i$ and $j$, so is the right
side.\qed

\vskip3ex

\noindent{\sl Proof of Lemma \ref{noname}.} We will prove this using the formula
\ref{torsionformula}.    We assume that the determinant of $I-\tilde{\gamma}$
(resp. $I-\tilde{\beta}$) is non-zero since if not    both sides of each
formula in Lemma \ref{noname} are zero from the definitions.

  We begin by
setting up the notation.  Let $F=\qq(t_1,\cdots,t_n)$.  From the definitions the Gassner and
reduced Gassner  representations are related by the diagram
\[
 \begin{diagram} 
\node{0}\arrow{e}\node{F}\arrow{e,t}{w}
       \node{H^1(X_0,p)}\arrow{e,t}{i^*}\arrow{s,l}{\gamma}
          \node{H^1(X_0)}\arrow{s,r}{\tilde{\gamma}}\arrow{e}\node{0}\\
\node{0}\arrow{e}\node{F}\arrow{e,t}{w}
       \node{H^1(X_0,p)}\arrow{e,t}{i^*} 
          \node{H^1(X_0)} \arrow{e}\node{0}
 \end{diagram}
\]
where the horizontal rows are exact,  $H^1(X_0,p)\cong F^n$, $H^1(X_0)\cong
F^{n-1}$. Fix the basis for $H^1(X_0,p)$ as in Section 4 and call it $\{e_0,\cdots,
e_n\}$.   and the map labeled $w$ takes $1\in F$ to the vector 
$$w=\pmatrix{1-t_1\cr\vdots\cr 1-t_n\cr}.$$
Recall from Proposition \ref{sadpoi} that  the
vector
$u=(t_1^{-1},\cdots , (t_1t_2\cdots t_n)^{-1})$ satisfies 
$u(I-\gamma)=0$ and that $(I-\gamma)w=0$.   Let $q:H^1(X_0,p)\to F$ denote
the map
$q(z_1,\cdots,z_n)=\sum_i u_i z_i$.

Denote the matrix   for $I-\gamma:H^1(X_0,p)\to
H^1(X_0,p)$ in the given basis by
$M$. 
 Clearly the image of $e_2,\cdots, e_n$ in $H^1(X_0)$ form a basis.   Use this basis
for $H^1(X_0)$ and continue to call it $e_2,\cdots,e_n$.  Then we have three based
acyclic chain  complexes and an exact sequence relating them described in the 
following diagram.
\begin{equation}\label{complexes}
 \begin{diagram}\dgARROWLENGTH=1.5ex
 \node{0}\arrow{e}\node{F}\arrow{e,t}{w_1}\arrow{s,l}{1}
\node{F}\arrow{s,l}{w_1^{-1}w}\arrow{e,t}{0}\node{F}\arrow{s,l}{w_1^{-1}w}\arrow{e,t}{w_1^{-1}q(w)}
\node{F}\arrow{e}\arrow{s,r}{1}\node{0}\\
\node{0}\arrow{e}\node{F}\arrow{s}\arrow{e,t}{w}
\node{H^1(X_0,p)}\arrow{e,t}{M}\arrow{s}\node{H^1(X_0,p)}\arrow{s}\arrow{e,t}{q}
\node{F}\arrow{s}\arrow{e}\node{0} \\
\node{0}\arrow{e}\node{0}\arrow{e}\node{H^1(X_0)}\arrow{e,t}{\widetilde{M}}
\node{H^1(X_0)}\arrow{e}\node{0}\arrow{e}\node{0}
\end{diagram}
\end{equation}
In Diagram \ref{complexes} the horizontal rows are based, acyclic complexes. Also,
$w_1=1-t_1$ is the first entry in the vector $w$.

For notational ease call the top complex $C_*$, the middle complex $D_*$, and
the bottom complex $E_*$.  The reader can easily check that:
\begin{enumerate}
\item the short exact sequence $0\to C_*\to D_*\to
E_*\to 0$ has compatible bases, and
\item the  torsion of $D_*$ is exactly the torsion computed in the proof of Lemma 
\ref{oinvxy}, so 
$$\tau(D_*)={{u_1 w_1}\over {\det(M(11))}}=\Delta_L^{-1}.$$
\end{enumerate}
 
Since the complexes are compatibly based and acyclic, the formula
\ref{torsionformula}  implies that
\begin{equation}\label{jeidk}
\tau(C_*)\tau(E_*)=\tau(D_*).
\end{equation}
But trivial computations show that  $\tau(E_*)= \det(\widetilde{M})^{-1}=
\det(I-\tilde{\gamma})^{-1}$ and 
$$\tau(C_*)=q(w)=\sum_i u_iw_i=\sum_i (t_1\cdots t_i)^{-1}(1-t_i)=
(t_1\cdots t_n)^{-1}-1.$$
Substituting into  Equation \ref{jeidk}  yields
$$((t_1\cdots t_n)^{-1}-1) \det(I-\tilde{\gamma})^{-1}=\Delta_L^{-1}.$$
This establishes the first statement in the Lemma.

For the second statement, repeat the argument, but this time with 
$$F=\qq(t), \ M=I-\beta,\ \widetilde{M}=I-\tilde{\beta}.
$$ Also set
$$ w=\pmatrix{1\cr\vdots\cr
1\cr}  \mbox{ and } u=(t^{-1}, t^{-2},\cdots, t^{-n})$$
since we saw in Section 6 that these are the left and right $1$ eigenvectors for the
Burau representation.

Then again $\tau(D_*)= \bar{\Delta}_L^{-1}$ but this time
$\tau(E_*)=\det(\widetilde{M})^{-1}=
\det(I-\tilde{\beta})^{-1}$ and
$$\tau(C_*)=q(w)=\sum_i u_iw_i=t^{-1}+t^{-2}+\cdots +t^{-n}.$$
Thus Equation \ref{jeidk} reads
$$(t^{-1}+t^{-2}+\cdots +t^{-n}) \det(I-\tilde{\beta})^{-1}=
\bar{\Delta}_L^{-1}$$
establishing the second claim. \qed

For the rest of this section $R=\Lambda=\zz[\zz^n]=\zz[t_i, t_i^{-1}]$ and
$F=\qq(t_1,\cdots, t_n)$  is its quotient field.

The following theorem is   ``folklore''; it is well known in the case of knots.  
An excellent exposition of Reidemeister torsion in knot theory can be found in
\cite{turaev}. There the connections between Alexander polynomials of links and
Reidemeister torsion are fully developed.

\begin{theorem}\label{lkwen}

 Let $\hat{L}\subset S^3$ be an $n$-component link; $n>1$.  Let
$C_*=C_*(S^3-\hat{L};\Lambda)$ be the free $\Lambda$ chain complex of  cellular chains
in the universal abelian (i.e. $\zz^n$) cover of $S^3-\hat{L}$.

 Then
$C_*\otimes F$ is acyclic  if and only if the Alexander polynomial
$\Delta_{\hat{L}}$ is non-zero, and in this case
$$\Delta_{\hat{L}}=\tau(S^3-\hat{L})^{-1}.$$

\end{theorem}

 Notice that when $\Delta_{\hat{L}}$ is zero, $\tau(S^3-\hat{L};K)^{-1}$
is a (non-zero) generalization of the Alexander polynomial which
depends in general on the choice of free submodule $K\subset
H_*(Z;\Lambda)$.  

\vskip3ex

\noindent{\sl Proof of Theorem \ref{lkwen}.}
We sketch the argument.  Let $Z$ denote
$S^3-\hat{L}$.  The Wirtinger presentation of $\pi_1(Z)$ gives a
2-complex cellular structure for $Z$ (up to homotopy) provided we
drop one of the relations, say the one coming from the first
crossing.  Take the basis ${\bf e}_i$ of $C_*$ to be the  cellular
basis.  But then the boundary map from $2$-chains to $1$-chains
in this basis is exactly the Alexander matrix with its first row
dropped, which we denote by $V'$.  The boundary map from $1$-chains to $0$-chains is
given by the $(c\times 1)$ matrix taking each 1-cell to the appropriate $1-t_i$.
Thus the chain complex $C_*=0\to C_2\to C_1\to C_0\to 0$ looks like

 \[
 \begin{diagram} 
\node{0}\arrow{e}\node{\Lambda^{c-1}}\arrow{e,b}{V'}\node{\Lambda^{c}}\arrow{e,b}
{({1-t_{\alpha(1)}},\cdots )^T}\node{\Lambda}\arrow{e}\node{0}
  \end{diagram}
\]

 Choosing   free submodules $K_i\subset H_i(Z;\Lambda)$ 
defines a torsion $\tau(Z;K)\in F^*/\Lambda^*$.

 Note that
$H_0(Z;F)=F/\langle 1-t_i\rangle=0$.   Also
$H_2(Z;F)$ is the kernel of $\partial_2$, and we have identified it
with the kernel of the Alexander matrix with its first row
dropped.  Recall that the Alexander polynomial $\Delta_{\hat{L}}$ 
is just $\det(V(ij))/(1-t_{\alpha(j)})$ (where the labeling of the
generators is taken so that the $j$th generator lies on the
$\alpha(j)$th component.)  Thus, $\Delta_{\hat{L}}\ne 0$ if and only
if
$\det(V(ij))\ne 0$ for each $i,j$  which is true if and only
$\partial_2$ is injective.  

Assume by relabeling if necessary that the first generator is part of
the first component for the string link, i.e. $\alpha(1)=1$.  
 Take  $\tilde{{\bf b}}_0$ to be the first element
$e_{1,1}$ of the basis ${\bf e}_1$, so that ${\bf
b}_0=(1-t_1)e_{1,1}$. Also take 
$\tilde{{\bf b}}_1 ={\bf e_2}$, so that ${\bf b}_1=V'{\bf e_2}$. 
Also

 Computing directly
with the definition one sees that if the Alexander polynomial is
non-zero, then
 \begin{eqnarray*}
\tau&=&{{[{\bf b}_0|{\bf e}_0]\ [{\bf \tilde{b}}_1|{\bf e}_2]}\over
{[{\bf b}_1|{\bf e}_1]}}\\
&=&{{(1-t_1)\cdot 1}\over{\det(V(11))}}\\
&=&\Delta_{\hat{L}}^{-1} 
\end{eqnarray*}
\qed

We next wish to give our generalization of Theorem \ref{factorization} using
torsion.  We first recall and set up some notation.

 Let $L\subset D^2\times I$ be an $n$-component
string link;
$n>1$.  Let $X$ denote the complement $(D^2\times I)-L$. Let $X_0$
and $X_1$ denote the punctured disks $( D^2-\{p_i\})\times \{0\}$
and $( D^2-\{p_i\})\times \{1\}$.  Let $W$ be the complement of
the closure  $\hat{L}$ viewed as a link in the solid torus, so  
$$W= X \cup X_0\times I$$
and let $Z$ be the complement of $\hat{L}$ in $S^3$, so that
$$Z=W\cup Q$$
where $Q$ is a solid torus.  

 Each of the spaces $X_0,X_1,X,W,Z$,
and $Q$ has a   $\Lambda$ (and hence also $F$)  coefficient system
defined by the abelianization maps on fundamental groups.   

\begin{lemma}  \label{cxvoihjtr}

\begin{enumerate}
\item The inclusion  $ W\subset Z$ induces  
isomorphisms  
$H_i(W;F)\to H_i(Z;F)$.  These groups are zero except possibly if
$i=1$ and $2$.
   Moreover,
\begin{equation}\label{poti}
\tau(Z,W)=(1-t_1\cdots t_n).
\end{equation}

\item The projections $X_0\times I\to  X_0$, $X_0\times I\to X_1$,
and the inclusion $X_0\subset X$ induce isomorphisms on
$F$-homology. The  homology groups $H_i(X_0;\Lambda)$ are zero except for
$i=1$ and $H_1(X_0;\Lambda)=\Lambda^{n-1}$. 
\end{enumerate}
\end{lemma}

\noindent{\sl Proof.} The excision theorem implies that  the  homology groups
$H_i(Z,W)$ are isomorphic to 
$H_i(Q,T)$, where
$T=\partial Q$ is a torus. The two generators of $\pi_1T$ can be
taken to be the curve $\mu$ bounding a disk in $Q$ and
$\lambda=\partial X_0$.   The representation  $\pi_1T\to \zz^n$ takes
$\mu$ to $1$  and $\lambda$ to $1-t_1\cdots t_n$.  Notice that $Q$
deforms to $\lambda$.   It follows readily that $H_i(T;F)$ and
$H_i(Q;F)$ vanish for all $i$. Thus $H_i(Q,T;F)=H_i(Z,W;F)=0$ for
all $i$ and $\tau(Z,W)=\tau(Q,T)$. At this point we leave to the
reader the elementary  computation 
$$
\tau(Z,W)=\tau(Q,T)=(1-t_1\cdots t_n).
$$

We have calculated in Section 2  that 
$$H_i(X_0;\Lambda)=\cases{\Lambda^{n-1}&if $i$=1\cr
         0 &otherwise.\cr}$$
and that with $F$ coefficients the inclusion $H_*(X_0;F)\to H_*(X;F)$
is an isomorphism.  The second assertion follows easily. \qed

The following notation will hold for the rest of this section.  Let  $K_i\subset
H_i(W;\Lambda), \ i=1,2$ be   free submodules so that the inclusion  induces   
isomorphisms  $K_i\otimes F\to H_i(W;F)$. For
convenience denote the image of
$K_i$ in  $H_i(Z;\Lambda)$ by $K_i$ also.  (Recall that by Lemma \ref{cxvoihjtr}
$K_i$ injects into  $H_i(Z;\Lambda)$.)  Let
$J=H_1(X_0\times I;\Lambda)$.  Also denote the (injective)
image of
$J$ in the $\Lambda$-homology groups of $X_0,X_1$ and $X$ by $J$.
 Notice that   $J$ is canonical, i.e. no choices are involved in making its definition. In
contrast, there is no natural choice of $K_i$ in general.

\begin{theorem}
  Let $\Phi$ denote Mayer-Vietoris diagram
\[
 \begin{diagram}\dgARROWLENGTH=1.6ex
 \node[2]{X_0\times I}\arrow{se}\\
\node{X_0\cup X_1}\arrow{ne}\arrow{se}\node[2]{W}\\
\node[2]{X}\arrow{ne}
 \end{diagram}
\]
 
 The
Mayer-Vietoris homology sequence with $F$ coefficients for $\Phi$ reduces to the
sequence 
\begin{equation}\label{ciuhnbg}
0\to H_2(W)\to H_1(X_0)\oplus H_1(X_1)\mapright{\xi}
H_1(X)\oplus H_1(X_0\times I)\to H_1(W) \to 0
\end{equation}
 Denote its torsion (which is defined once one specifies the spaces $K_i$ and $J$ in the
manner described above)  by
$\tau({\cal H}(\Phi);K)$.  

Then the map $\xi$ in the sequence (\ref{ciuhnbg}) is given by the matrix
$$\xi=\pmatrix{I&I\cr \tilde{\gamma}&I\cr}$$
 where $\tilde{\gamma}$ denotes the reduced Gassner representation.
 
Furthermore,
\begin{equation}\label{oikjvnvcx}
\tau(Z;K)^{-1}={{\tau(L)\ \tau({\cal H}(\Phi);K)}\over {  1-t_1\cdots t_n   }}
\end{equation}
Where $\tau(L)$ is the torsion of the string link (defined in Section 6).

 Finally, if $\Delta_L\ne 0$  then $K=0$,
$${{\tau(H(\Phi);K)}/{1-t_1\cdots t_n}}=\Delta_L,$$
$$\tau(Z;K)^{-1}=\Delta_{\hat{L}},$$  and Equation
\ref{oikjvnvcx}
 reduces to Theorem \ref{factorization}.

\end{theorem}

\noindent{\sl Proof.}
Lemma \ref{cxvoihjtr}  shows that all the other groups in the Mayer-Vietoris
sequence for $\Phi$ vanish.   
From the definition of the reduced (homology) Gassner matrix the map labeled $\xi$ in
the sequence
\ref{ciuhnbg} is given in this basis by the matrix
$$\xi=\pmatrix{I&I\cr \tilde{\gamma}&I\cr}.$$

Applying Equation \ref{torsionformula} to the Mayer-Vietoris
sequence and using the fact that the torsion of a direct sum is the
product of torsions, one concludes that
\begin{equation}\label{hvn}
\tau({\cal H}(\Phi);K)={{\tau(X_0\times I;J)\ \tau(X;J)}\over{\tau(X_0;J)\
\tau(X_1;J)\ \tau(W;K)}}.
\end{equation}

Notice, however, that $X_0\times I$ deforms to $X_0$ and $X_1$, and
by the way we chose the bases for $J$ it follows that
$\tau(X_0\times I;J)=\tau(X_0;J)=\tau(X_1;J)$. Thus Equation
\ref{hvn} reduces to 

\begin{equation}\label{hvni}
\tau({\cal H}(\Phi);K)={{\tau(X;J)}\over {\tau(X_1;J)\ \tau(W;K)}}.
\end{equation}

Now consider the pair $(Z,W)$. The complex
$C_*(Z,W;F)$ is acyclic and  the choice of bases for $K_i\subset
H_i(W;\Lambda)$ were used to define the bases for $H_i(W;F)$ {\it and}
$H_i(Z;F)$, it follows that the torsion of the long exact sequence
in homology for the pair $(Z,W)$ is equal to $1$; in fact the
sequence reduces to isomorphisms which are the identity in the given
bases.  Thus Equation \ref{torsionformula} implies that
$$1={{\tau(W;K)\ \tau(Z,W)}\over {\tau(Z;K)}}$$
and so 
\begin{equation}\label{hvnii}
\tau(W;K)= \tau(Z;K) \tau(Z,W)^{-1}.
\end{equation}

Finally consider the pair $(X,X_1)$.  A similar argument as the
previous paragraph shows that
\begin{equation}\label{hvniii}
1={{\tau(X_1;J)\ \tau(X,X_1)}\over{\tau(X;J)}}.
\end{equation}

 Combining Equations \ref{poti},  \ref{hvni}, \ref{hvnii},  and
\ref{hvniii} yields 
\begin{equation}\label{hvniv}
\tau(Z;K)^{-1}={{\tau(X,X_1)^{-1}\ \tau({\cal H}(\Phi);K)}\over {  1-t_1\cdots
t_n   }}.
\end{equation}

 Moreover
$\tau(X,X_1)=\tau(L)^{-1}$ 

Since we defined
$\tau(L)$ using the (dual) cochain complex.    Substituting this   into Equation
\ref{hvniv} establishes Equation \ref{oikjvnvcx}.

Now recall that $\Delta_L$ is non-zero if and only if
$\det(I-\tilde{\gamma})$ is non-zero.  Since
$$\det(I-\tilde{\gamma})=\det\pmatrix{I&I\cr \tilde{\gamma}&I\cr},$$
it follows that $\Delta_L$ is non-zero if and only if
$H_*(Z;F)=H_*(W;F)=K=0$ and in this case the sequence \ref{ciuhnbg}
shows that 
$$\tau({\cal H}(\Phi))= \det(I-\tilde{\gamma}).  $$ Thus
Equation \ref{hvniv} reduces in this case to
\begin{equation}\label{hvnv}
\tau(Z)^{-1}={{\tau(L)\ \det(I-\tilde{\gamma})}\over {  1-t_1\cdots t_n   }}
\end{equation}

Lemma \ref{noname} shows that 
$\Delta_L={{\det(I-\tilde{\gamma})}\over{((t_1\cdots t_n)^{-1}-1)}}$.
Notice that $$1-t_1\cdots t_n =(t_1\cdots t_n)((t_1\cdots t_n)^{-1}-1).$$ 
Theorem \ref{lkwen} implies that $\tau(Z )^{-1}=\Delta_{\hat{L}}$. 

Substituting these into  Equation
\ref{hvnv}  and working modulo units in $\Lambda$ one obtains, when $\Delta_L\ne 0$, the
formula
$$\Delta_{\hat{L}}= \tau(L)\ \Delta_L.$$
This is just 
 the formula of Theorem
\ref{factorization}.

\qed

 \vskip3ex 

\noindent Department of Mathematics \hfill pkirk@indiana.edu

\noindent Indiana University\hfill livingst@indiana.edu

\noindent Bloomington, IN 47405 \hfill zhewang@indiana.edu


\begin{thebibliography}{99}

\bibitem{abdul}M. N. Abdulrahim, {\em Complex specializations of the reduced
Gassner representation of the pure braid group,}
Proc. Amer. Math. Soc. 125 (1997), no. 6, 1617--1624. 

\bibitem{aps2} M.F. Atiyah, V.K. Patodi, and I.M. Singer {\em Spectral asymmetry and
Riemannian geometry. II,} Math Proc. Camb. Phil. Soc. 78 (1975), 405--432.


\bibitem{barnatan} D. Bar-Natan {\em On the Vassiliev knot invariants,}
Topology 34 (1995), no. 2, 423--472. 

\bibitem{birman} J. Birman  {\em Braids, links, and mapping class groups,}
Annals of Mathematics Studies No. 82, Princeton University Press, Princeton,
N.J., 1974.

\bibitem{orr} T. Cochran and K. Orr, {\em Not all links are concordant to
boundary links,} Ann. of Math. (2) 138 (1993), no. 3, 519--554. 

\bibitem{cochran} T. Cochran,
{\em Non-trivial links and plats with trivial Gassner matrices,}
Math. Proc. Cambridge Philos. Soc. 119 (1996), no. 1, 43--53. 


\bibitem{fox2} R. Fox  {\em Free differential calculus. II. The
isomorphism problem of groups,} Ann. of Math. (2) 59, (1954) 196--210. 

\bibitem{fox} R. Fox {\em A quick trip through knot theory,} 1962 Topology of
3-manifolds and related topics (Proc. The Univ. of Georgia Institute, 1961)
 120--167 Prentice-Hall, Englewood Cliffs, N.J.

\bibitem{fox-milnor} R.  Fox and J. Milnor {\em Singularities of $2$-spheres
in
$4$-space and cobordism of knots,} Osaka J. Math. 3 (1966) 257--267.

\bibitem{habegger} N. Habegger and X. S. Lin {\em The classification of links
up to link-homotopy,} J. Amer. Math. Soc. 3 (1990), no. 2, 389--419.

\bibitem{jones} V. Jones, {\em A polynomial invariant for knots
via von Neumann algebras,} Bull. Amer. Math. Soc. 12 (1985), no. 1,
103--111.  

\bibitem{kirk-livingston} P. Kirk and C. Livingston {\em Twisted Alexander
invariants, Reidemeister torsion, and the Casson-Gordon invariants,} to appear,
Topology (1996).

\bibitem{ledimet} J.-Y.  Le Dimet
{ \em Enlacements d'intervalles et representation de Gassner,}  
Comment. Math. Helv. 67 (1992), no. 2, 30-6-315.

\bibitem{levine} J. Levine {\em A factorization of the Conway polynomial,}
preprint (1997). 

\bibitem{lin-wang} X.S. Lin, F. Tian, and Z. Wang {\em Burau representations
and random walks on string links,} Pacific J. of Math 182 (1998) 289--302.

\bibitem{long} Long, D. D. {\em On the linear representation of braid groups,}
Trans. Amer. Math. Soc. 311 (1989), no. 2, 535--560.

\bibitem{milnor} J. Milnor {\em Whitehead torsion,} Bull. Amer. Math. Soc. 72
1966 358--426.

\bibitem{moody}  J. Moody {\em
The faithfulness question for the Burau representation,} 
Proc. Amer. Math. Soc. 119 (1993), no. 2, 671--679.

\bibitem{penne}  R. Penne {\em Multi-variable Burau matrices and labeled
line configurations,}  J. Knot Theory Ramifications
4 (1995), no. 2, 235--262.

\bibitem{stallings}  J. Stallings {\em Homology and central series of groups,}
J. Algebra 2 (1965) 170--181.

\bibitem{squier} C. Squier,{\em The Burau representation is unitary,} Proc.
Amer. Math. Soc. 90 (1984), no. 2, 199--202.

\bibitem{turaev} V. Turaev , {\em Reidemeister torsion in knot theory,} Russian
Math. Surveys 41:1 (1986), 119--182.

\end{thebibliography}
 \end{document}